\documentclass{amsart}

\begin{document}
\title{Regular Fredholm Pairs}
\author{Enrico Boasso}
\begin{abstract}In this work it is introduced the notion of regular Fredholm pair, i.e. a 
Fredholm 
pair whose operators are regular. The main properties of these
objects are studied, and what is more, they are entirely classified. Furthermore, 
the index of a Fredholm pair turns out to be an extremely useful tool in the 
description of
the aforementioned objects. Finally, regular Fredholm pairs are characterized in terms 
of regular Fredholm symmetrical pairs, exact chains of multiplication operators, 
and invertible  Banach space operators. \par
\vskip.5cm
\noindent KEYWORDS: Index, Fredholm pairs, regular operators.\par
\noindent MSC (2000): Primary 47A13; Secondary 47A53.
\end{abstract}
\maketitle

\noindent \bf{1. Introduction}\rm \vskip.3cm
\indent There are many ways to extend the notions of Fredholm operator and
index to several variable operator theory. For instance, both Fredholm complexes of
Banach spaces and the related notion of Fredholm pair have associated
an index with good stability properties, see for example \cite{9}, \cite{1}, \cite{2}, 
\cite{3}. \par

\indent On the other hand, regular maps are a natural generalization of Fredholm 
ope-
rators. However, the boundary maps of a Fredholm complex of Banach spaces are 
generally not regular. In fact, in order for such a complex
to have this property, it must be a split Fredholm complex, see for example
\cite{4}. As regard Fredholm pairs, since the operators of such a pair are generally 
not regular,
a similar situation is encountered. The main objective of this work consists
in the study of regular
Fredholm pairs, i.e. Fredholm pairs whose operators are regular.\par

\indent In the next section the notion of regular Fredholm pair is introduced. Moreover, 
some definitions and facts needed for the present work are reviewed, and 
some preliminary 
and general results regarding regularity are also proved. In section three the 
objects under
consideration 
are entirely classified. In section four the index turns out to be an 
extremely useful tool to  describe regular Fredholm pairs. Weyl pairs 
are also introduced and considered.      
Finally, in section five three
characterizations of regular Fredholm pairs are proved. In fact, 
these objects are characterized in terms of regular symmetrical Fredholm pairs, 
exacts chains of
multiplication operators, and invertible 
Banach space operators.   
\par           

\indent The author wishes to express his indebtedness to Professors C.-G. 
Ambrozie and
F.-H. Vasilescu. These researchers sent kindly to the author several works 
authored by them, 
which have been necessary for the elaboration of the present article.  \par
\newpage
\noindent \bf{2. Notations, Definitions and Preliminary Results}\rm\vskip.3cm

\indent From now on $X$ and $Y$ denote two Banach spaces, $L(X,Y)$ the algebra of 
all linear and continuous operators defined on $X$ with values in $Y$, and $K(X,Y)$ 
the closed
ideal of all compact operatots of $L(X,Y)$. As usual, when $X=Y$, $L(X,X)$ and 
$K(X,X)$ are denoted
by $L(X)$ and $K(X)$ respectively. For
every $S\in L(X,Y)$ the null space of $S$ is denoted by
$N(S)=\{x\in X\colon \hbox{  }S(x)=0\}$, and the range of $S$ by $R(S)=
\{ y\in Y: \hbox{  }\exists
\hbox{ } x\in X \hbox{ such that }y=S(x)\}$. Next follows the definition of Fredholm pair, 
see for instance \cite{2}.\par

\newtheorem*{def2.1}{Definition 2.1} 
\begin{def2.1} Let $X$ and $Y$ be two Banach spaces and let $S\in L(X,Y)$,
$T\in L(Y,X)$ be such that the following dimensions are finite:
\begin{align*}
&a\colon =\dim N(S)/(N(S)\cap R(T)),& &b\colon = \dim R(T)/(N(S)\cap R(T)),\\
&c\colon = \dim N(T)/(N(T)\cap R(S)),& &d\colon =\dim R(S)/(N(T)\cap R(S)).
\end{align*}
A pair $(S,T)$ with the above properties is called a Fredholm pair.\par
\indent Let $P(X,Y)$ denote the set of all Fredholm pairs. If $(S,T)\in P(X,Y)$, then
the index of $(S,T)$ is defined by the equality
$$
\hbox{\rm ind} \hbox{ } (S,T)\colon = a -b -c +d.
$$
\indent In particular, if $(S,T)\in P(X,Y)$ is such that $ST=0$
and $TS=0$, that is if $b=d=0$, then $(S,T)$ is said
a Fredholm symmetrical pair, see \cite{10}. Note that in this case
$(S,T)$ and $(T,S)$ are Fredholm chains, see \cite[10.6]{5} and \cite{6}.\par 
\end{def2.1}
\indent Before going on, several properties of Fredholm pairs are recalled, see \cite{2}.
\par
\newtheorem*{rem2.2}{Remark 2.2}
\begin{rem2.2}\rm
First of all, observe that if $S\in L(X,Y)$ is a Fredholm operator, then 
$(S,0)$ is a Fredholm pair. Furthermore, $\hbox{\rm ind}\hbox{ } S=\hbox{\rm ind}
\hbox{ } (S,0)$. Consequently, the definition of Fredholm pair extends the notion of 
Fredholm operator to several variable operator theory.\par
\indent In second place, note that if $(S,T)\in P(X,Y)$, then $(T,S)\in P(Y,X)$ and
$$
\hbox{\rm ind}\hbox{ } (T,S)=-\hbox{\rm ind}\hbox{ } (S,T).
$$ 
\indent In third place, if $(S,T)\in P(X,Y)$, then $R(S)$ and $N(T) +R(S)$ are closed 
subspaces in $Y$. Similarly, $R(T)$ and $N(S)+R(T)$ are closed subspaces in $X$.
What is more, there are finite dimensional subspaces in $X$, $X_1$ and $X_2$,
and in $Y$, $Y_1$ and $Y_2$, such that:\par
\it{i)}\rm{ }$\dim X_1=a$, $\dim X_2=b$, $\dim Y_1=c$, $\dim Y_2=d$,\par

\it{ii)}\rm{ }$N(S)= (N(S)\cap R(T))\oplus X_1$, $R(T)=(N(S)\cap R(T))\oplus X_2$,\par

\it{iii)}\rm{ }$N(T)= (N(T)\cap R(S))\oplus Y_1$, $R(S)=(N(T)\cap R(S))\oplus Y_2$.
 
\noindent In particular,
$$
N(S) + R(T)=(N(S)\cap R(T))\oplus (X_1\oplus X_2),
$$
and
$$
N(T) + R(S)=(N(T)\cap R(S))\oplus (Y_1\oplus Y_2).
$$
\indent Moreover, $S$  induces an isomorphism
$$
X_2= R(T)/(N(S)\cap R(T))\xrightarrow{\cong} R(ST).
$$
Similarly, $T$ induces and isomorphism
$$
Y_2= R(S)/(N(T)\cap R(S))\xrightarrow{\cong} R(TS).
$$
In particular, $\dim R(ST)=b$ and $\dim R(TS)=d$.\par

\indent On the other hand, it is easy to prove that if $\tilde{Y}_2$ is another subspace
such that $R(S)=(N(T)\cap R(S))\oplus \tilde{Y}_2$, then
$\dim \tilde{Y}_2=d$ and $T$ induces an isomorphism
$$
\tilde{Y}_2= R(S)/(N(T)\cap R(S))\xrightarrow{\cong} R(TS).
$$ 
\indent Finally, interchanging $S$ and $T$, similar properties for the operator $S$ 
can be proved.\par
\end{rem2.2}
\indent Next follows the definition of regular operator, see 
for example \cite{5}.\par

\newtheorem*{def2.3}{Definition 2.3} 
\begin{def2.3}
Let $X$ and $Y$ be two Banach spaces 
and let $T\in L(X,Y)$. The operator $T$ is called regular or relatively
Fredholm, if there is $T' \in L(Y,X)$ for which
$$
T=TT'T.
$$
\end{def2.3} 
\indent If $T$ is a regular bounded and linear map, the operator
$T'$ in Definition 2.3 is called a \it{generalized inverse}\rm, or
\it{pseudo inverse}\rm, for $T$.  If, in addition, $T$ is a generalized
inverse for $T'$, that is if
$$
T'=T'TT',
$$
then $T'$ is called a \it{normalized} \rm generalized inverse,
see \cite[3.8]{5} and \cite{6}. It is well known that if $T'$ is a generalized
inverse for $T$, then 
$$
T''=T'TT'
$$
is a normalized generalized inverse for $T$, see \cite[3.8]{5} and \cite{6}\par

\indent On the other hand, when the range of $T$ is closed,
the condition of being a regular operator is equivalent to the fact that $N(T)$ and
$R(T)$ are complemented subspaces in $X$ and $Y$
respectively, see \cite[3.8.2]{5}.\par

\indent In the following proposition Fredholm pairs
whose operators have complemented ranges and
null spaces are studied.\par

\newtheorem*{prop2.4}{Proposition 2.4} 
\begin{prop2.4}Let $X$ and $Y$ be two Banach spaces 
and $(S,T)\in P(X,Y)$. The following assertions are
equivalent:\par
 i) $R(T)$ is a complemented subspace in $X$;\par
ii) $N(S)$ is a complemented subspace in $X$;\par
iii) $N(S)+R(T)$ is a complemented subspace in $X$;\par
 iv) $N(S)\cap R(T)$ is a complemented subspace in $X$.\par 
\indent Similarly, the following assertions are equivalent:\par
i) $R(S)$ is a complemented subspace in $Y$;\par
ii) $N(T)$ is a complemented subspace in $Y$;\par
iii) $N(T)+R(S)$ is a complemented subspace in $Y$;\par
 iv) $N(T)\cap R(S)$ is a complemented subspace in $Y$.\par 
\indent In particular, if $S$ (resp. $T$) is a regular operator,
then $T$ (resp. $S$) also is a regular operator.\par
\end{prop2.4}
\begin{proof}
\indent According to Remark 2.2, all the subspaces involved
in the first part of the proposition are closed. Moreover, there are finite
dimensional subspaces $X_1$ and $X_2$ such that
$$
N(S)= (N(S)\cap R(T))\oplus X_1,\hskip.5cm R(T)=(N(S)\cap R(T))\oplus X_2,
$$ 
and 
$$
N(S) + R(T)=(N(S)\cap R(T))\oplus (X_1\oplus X_2).
$$
\indent Now well, it is clear that $i)$, $ii)$, and $iii)$ imply $iv)$.
On the other hand, according to \cite[6.3.5]{5}, $iv)$ implies $i)$, $ii)$ and $iii)$.\par
\indent A similar argument proves 
the second part of the proposition.\par
\end{proof}

\indent Next follows the definition of regular Fredholm pair.\par

\newtheorem*{def2.5}{Definition 2.5} 
\begin{def2.5} Let $X$ and $Y$ be two Banach spaces and let $S\in L(X,Y)$,
$T\in L(Y,X)$ be such that $(S,T)\in P(X,Y)$. If the operators $S$ and $T$ are regular, 
that is if $S$ and $T$ have the equivalent
properties of Proposition 2.4, then $(S,T)$ is called a regular Fredholm pair.\par
\indent The set of all regular Fredholm pairs is denoted
by $RP(X,Y)$.\par
\indent In particular, if $(S,T)\in RP(X,Y)$ is a Fredhlom symmetrical pair,
then $(S,T)$ is said a regular Fredholm symmetrical pair. Note that in this case
$(S,T)$ and $(T,S)$ are regular Fredholm chains, see \cite[10.6]{5} and \cite{6}.
\end{def2.5}

\newtheorem*{rem2.6}{Remark 2.6}
\begin{rem2.6}\rm  Note that, according to Proposition 2.4, 
if $(S,T)\in P(X, Y)$, then the property of being a regular Fredholm pair is equivalent to 
the fact that either the operator $S$ or the operator $T$ is regular. Furthermore, if 
$(S,T)\in P(X,Y)$, 
then there are closed subspaces $\tilde{X}$ and $\tilde{Y}$, in
$X$ and $Y$ respectively, such that 
\begin{align*}
X&=(R(T)\oplus X_1)\oplus \tilde{X}=(N(S)\oplus X_2)\oplus \tilde{X}=(N(S)+R(T))\oplus 
\tilde {X}\\
&=((N(S)\cap R(T))\oplus (X_1\oplus X_2))\oplus \tilde{X},
\end{align*}
and
\begin{align*}
Y&=(R(S)\oplus Y_1)\oplus \tilde{Y}=(N(T)\oplus Y_2)\oplus \tilde{Y}=(N(T)+R(S))\oplus 
\tilde {Y}\\
&=((N(T)\cap R(S))\oplus (Y_1\oplus Y_2))\oplus \tilde{Y},
\end{align*}
where $X_j$ and $Y_j$, $j=1$, $2$, are the finite dimensional subspaces
considered in Remark 2.2.\par

\indent In addition, if $S_2$ and $\tilde{S}$ are the restrictions of $S$ to
$X_2$ and to $\tilde{X}$ respectively, and if $\mathcal{S}=S_2\oplus \tilde{S}$, 
then according to Remark 2.2,
$$
\mathcal{S}\colon X_2\oplus \tilde{X}\xrightarrow{\cong} R(S),\hbox{ } 
R(S)=S(X_2)\oplus S(\tilde{X})=R(ST)\oplus R(\tilde{S}).
$$
\indent Similarly, if $T_2$ and $\tilde{T}$ are the restrictions of $T$ to
$Y_2$ and to $\tilde{Y}$ respectively, and if $\mathcal{T}= T_2\oplus \tilde{T}$, 
then according to Remark 2.2,
$$
\mathcal{T}\colon Y_2\oplus \tilde{Y}\xrightarrow{\cong} R(T),\hbox{ } R(T)=
T(Y_2)\oplus T(\tilde{Y})=R(TS)\oplus R(\tilde{T}).
$$
\end{rem2.6}
\newtheorem*{rem2.7}{Remark 2.7}
\begin{rem2.7}\rm Four examples of regular Fredholm pairs will be
considered. In first place, let $X$ and $Y$ be two Banach spaces and $(S,T)$ 
belong to $P(X,Y)$. According to \cite{2}, or Remark 2.2, and to \cite[3.8.2]{5},
it is clear that if $X$ and $Y$ are Hilbert spaces,
then $P(X,Y)=RP(X,Y)$. On the other hand, according to \cite[6.3.4]{5}, if $R(S)$
and $R(T)$ are finite dimensional subspaces of the Banach spaces $Y$ and 
$X$ respectively, then $(S,T)$ belongs to $RP(X,Y)$.
\end{rem2.7}

\indent Next consider $(\mathcal{X},d)$ a \it complex of Banach spaces \rm of finite 
length, that is a sequence 
$$
0\to X_n\xrightarrow{d_n} X_{n-1}\to\ldots\to X_1\xrightarrow{d_1}X_0\to 0,
$$
where $X_p$ are Banach spaces, $p=0,\ldots ,n$, and the bounded operators $d_p
\in L(X_p,X_{p-1})$ are such that $d_{p-1}d_p=0$, $p=1,\ldots ,n$. 
Define the homology groups of $(\mathcal{X},d)$ as $H_p(\mathcal{X},d)= 
N(d_p)/R(d_{p+1})$,
$p=0,\ldots ,n$. A complex $(\mathcal{X},d)$ is said \it Fredholm  \rm if all its homology 
groups are finite dimensional. If this is the case, then it is possible to associate 
to $(\mathcal{X},d)$ the 
integer
$$
\hbox{\rm{ind} } (\mathcal{X},d)=\sum_{p=0}^n (-1)^p \dim H_p (\mathcal{X},d),
$$    
which is called the \it{index }\rm or the \it{Euler characteristic }\rm of $(\mathcal{X},d)$.
\par

\indent In \cite{3} a Fredholm symmetrical pair was associated to each Fredholm complex.
In fact, as above consider a complex of Banach spaces $(\mathcal{X},d)$ and set
$$
X=\bigoplus_{p=2k}X_p,\hskip1cm Y=\bigoplus_{p=2k+1}X_p,
$$
and define the maps $S\in L(X,Y)$ and $T\in L(Y,X)$ as follows:
$$
S=\bigoplus_{p=2k}\hbox{ } d_p,\hskip1cm T=\bigoplus_{p=2k+1}\hbox{ } d_p,
$$
where $k\ge 0$, and $X_p=0$ when $p\ge n+1$.\par

\indent Since $(\mathcal{X},d)$ is a complex, $TS=0$ and $ST=0$. Furthermore, 
$(\mathcal{X},d)$
is a Fredholm complex if and only if $\dim$ $N(S)/R(T)$ and $\dim$ $N(T)/R(S)$
are finite dimensional, which is equivalent to the fact that $(S,T)$ is a Fredholm 
symmetrical pair. In addition,
$$
\hbox{\rm ind } (\mathcal{X},d)=\hbox{\rm ind } (S,T).
$$
\indent Now well, a complex of Banach spaces $(\mathcal{X},d)$ is called 
\it Fredholm split \rm
 if there are continuous linear operators
$$
X_{p-1}\xrightarrow{h_{p-1}}X_p\xrightarrow{h_p}X_{p+1},
$$
such that 
$$
d_{p+1}h_p+h_{p-1}d_{p}=I_p-k_p,
$$
where $k_p\in K(X_p)$, $p=0,\ldots ,n$. When $k_p= 0$ for $p=0,\ldots , n$,
$(\mathcal{X}, d)$ is said a split complex. According to \cite[2.7]{4}, it is easy to prove that
a complex $(\mathcal{X},d)$ is Fredholm split if and only if the above associated 
pair $(S,T)$ is a 
regular Fredholm symmetrical pair. \par

\indent Finally, consider $(\mathcal{K},\delta)$ a \it Fredholm chain of Banach 
spaces\rm, that is a sequence 
$$
0\to \mathcal{K}_n\xrightarrow{\delta_n} \mathcal{K}_{n-1}\to\ldots\to \mathcal{K}_1
\xrightarrow{\delta_1} \mathcal{K}_0\to 0,
$$
where $\mathcal{K}_p$ are Banach spaces, and the bounded operators
$\delta_p\in L(\mathcal{K}_p,\mathcal{K}_{p-1})$ are such that
$$
N(\delta_p)/(N(\delta_p)\cap R(\delta_{p+1}))\hbox{ and }
R(\delta_{p+1})/(N(\delta_p)\cap R(\delta_{p+1}))
$$   
are finite dimensional subspaces of $\mathcal{K}_p$, $p=0,\ldots ,n$.\par 
\indent Recall that in \cite{7} it was introduced the more general concept of 
\it semi-Fredholm chains\rm.
However, since the main concern of this article consists in Fredholm objects,
only Fredholm chains will be considered. Furthermore,
observe that since
$\dim R(\delta_{p-1}\delta_p)$ is finite dimensional, $p=1,\ldots ,n$, a Fredholm
chain $(\mathcal{K},\delta)$
is a particular case of what in \cite{8} was called an \it essential complex of Banach spaces. \rm   
 \par
\indent As in the case of Fredholm complexes of Banach spaces, it is possible
to associate an index to any Fredholm chain. In fact, if $(\mathcal{K},\delta)$ 
is such an object, then define its index as
\begin{align*}
\hbox{ind } (\mathcal{K},\delta)=\sum_{p=0}^n&(-1)^p (\dim\hbox{ }
N(\delta_p)/(N(\delta_p)\cap R(\delta_{p+1}))\\
&-\dim R(\delta_{p+1})/(N(\delta_p)\cap R(\delta_{p+1}))),
\end{align*}
see \cite{7}.\par
\indent Now well, given a Fredholm chain $(\mathcal{K},\delta)$, define $X$, $Y$,
$S\in L(X,Y)$ and $T\in L(Y,X)$ as it has been done above for a Fredholm
complex of Banach spaces. Then, as in \cite{7}, 
$(\mathcal{K},\delta)$ is a Fredholm chain if and only if the associated pair $(S,T)$ is a
Fredholm pair. In addition, 
$$
\hbox{ind }(\mathcal{K},\delta) =\hbox{ind } (S,T).
$$

\indent In this work, in order to keep the analogy with complexes of Banach spaces, 
it will be said that a Fredholm chain $(\mathcal{K},\delta)$ is \it split \rm if
$N(\delta_p)$ is a complemented subspaces of $\mathcal{K}_p$,
$p=0,\ldots ,n$. Note that, as  in Proposition 2.4, this condition is equivalent to
the fact that $R(\delta_{p+1})$, or $N(\delta_p) + R(\delta_{p+1})$, or 
$N(\delta_p)\cap R(\delta_{p+1})$ is a complemented subspace of $\mathcal{K}_p$, 
$p=0,\ldots ,n$. Moreover, thanks to \cite[2.3]{8}, a split Fredholm
chain is a \it Fredholm essential complex \rm in the sense of \cite[2.2]{8}.\par

\indent Now well, it is not difficult to prove that $(\mathcal{K},\delta)$ is a split Fredholm
chain if and only if the above associated pair $(S,T)$ is a regular Fredholm pair.\par        

\newtheorem*{rem2.8}{Remark 2.8}
\begin{rem2.8}\rm Given Banach
spaces $X$ and $Y$, and $T\in L(X,Y)$ a Fredholm ope-
rator, then \it{any} \rm pseudo 
inverse $T'$ for $T$ is a Fredholm operator 
and ind $(T')=-$ ind $(T)$, see \cite[6.4.4]{5} and \cite[6.5.5]{5}. Nevertheless, 
as the following example
shows, these results do not hold any more for regular 
Fredholm pairs. In fact,
there is a regular Fredholm pair $(S,T)$, with $S'\in L(Y,X)$ and $T'\in L(X,Y)$ pseudo 
inverses for
$S$ and $T$ respectively, such that $(S',T')$ does not belong to $P(Y,X)$. 
\end{rem2.8}

\indent Consider finite dimensional Banach spaces $X_j$ and $Y_j$,
$j=1$, $2$, such that $\dim X_2= \dim Y_2$, and Banach spaces  
$\tilde{X}$, $\tilde{Y}$ and $N_j$, $j=1$, $2$, such that there are isomorphic operators
$$
\tilde{S}\colon \tilde{X}\xrightarrow{\cong} N_2,\hskip1cm \tilde{T}\colon\tilde{Y}
\xrightarrow{\cong}N_1.
$$
For example, take $\tilde{X}=N_2$ and $\tilde{Y}=N_1$, and $S$ and $T$ the
identity map of $\tilde{X}$ and $\tilde{Y}$ respectively.\par
\indent Define the Banach spaces
$$
X=X_1\oplus N_1\oplus X_2\oplus \tilde{X},\hskip1cm Y=Y_1\oplus N_2\oplus 
Y_2\oplus \tilde{Y},
$$
and the linear continuous maps $S\in L(X,Y)$ and $T\in L(Y,X)$ as follows:
\begin{align*}
&S\mid_{X_1\oplus N_1}\equiv 0,& &S\mid_{X_2}=S_2\colon X_2\to Y_2,&
&S\mid_{\tilde{X} }=\tilde{S}\colon \tilde{X}\to N_2, \\
&T\mid_{Y_1\oplus N_2}\equiv 0,& &T\mid_{Y_2}=T_2\colon Y_2\to X_2,&
&T\mid_{\tilde{Y} }=\tilde{T}\colon \tilde{Y}\to N_1, 
\end{align*}
where $S_2$ and $T_2$ are any isomorphic maps.\par
\indent It is easy to prove that $(S,T)\in RP(X,Y)$ and that ind $(S,T)=
\dim X_1-\dim Y_1$.\par
\indent On the other hand, consider the following operators $T'\in L(X,Y)$ and
$S' \in L(Y,X)$:
\begin{align*}
&T'\mid_{X_1\oplus\tilde{X}}\equiv 0,&
&T'\mid_{N_1\oplus X_2}=(T_2\oplus \tilde{T})^{-1}\colon N_1\oplus X_2\to
Y_2\oplus \tilde{Y},\\
&S'\mid_{Y_1}\equiv 0,&
&S'\mid_{N_2\oplus Y_2}=(S_2\oplus \tilde{S})^{-1}\colon N_2\oplus Y_2\to
X_2\oplus \tilde{X},
\end{align*}
and $S'\mid_{\tilde Y}\colon \tilde{Y}\to N_1$ any isomorphism.\par
\indent An easy calculation proves that $T'$ is a normalized generalized
inverse for $T$ and $S'$ is a pseudo inverse for $S$. However, since
$$
R(S')=N_1\oplus X_2\oplus\tilde{X}, \hbox{ } N(T')=X_1\oplus \tilde{X},\hbox{ }
N(T')\cap R(S')= \tilde{X},
$$
it is clear that
$$
R(S')/(N(T')\cap R(S'))=X_2\oplus N_1.
$$
\indent Therefore, if $N_1$, which is isomorphic to $\tilde{Y}$, is an infinite 
dimensional Banach space, then $(S',T')$ does not belong to $P(Y,X)$.\par
  
\indent Nevertheless, in the following proposition it is proved that 
if $(S,T)\in RP(X,Y)$, then there always exist normalized generalized 
inverses for $S$ and $T$, $S'$ and $T'$ respectively, such that
$(S',T')\in RP(Y,X)$.\par

\newtheorem*{prop2.9}{ Proposition 2.9}
\begin{prop2.9} Let $X$ and $Y$ be two Banach spaces and let 
$(S,T)\in RP(X,Y)$. Then, there is $(S',T')\in RP(Y,X)$ such that:\par
i) {\rm ind} $(S',T')=-$ {\rm ind} $(S,T)$;\par
ii) $S'$ and $T'$ are normalized generalized inverses for $S$ and
$T$ respectively.
\end{prop2.9}
\begin{proof}
   
\indent Consider $X$, $Y$, $S$ and $T$ presented as in Remark 2.6, and
define the operators $S'$ and $T'$ as follows:
\begin{align*}
&S'\mid_{ Y_1\oplus\tilde{Y}}\equiv 0,& &S' \mid_{ R(S)}=\iota_1\circ 
{\mathcal{ S}}^{-1}\colon R(S)\to X,\\
&T'\mid_{ X_1\oplus\tilde{X}}\equiv 0,& &T' \mid_{ R(T)}=\iota_2\circ 
{\mathcal{ T}}^{-1}\colon R(T)\to Y,
\end{align*}
where $\iota_1\colon X_2\oplus\tilde{X}\to X$ and $\iota_2\colon Y_2\oplus \tilde{Y}
\to Y$ are the natural inclusion maps. \par

\indent A straightforward calculation proves that $S'$ and $T'$ are
normalized generalized inverses for $S$ and $T$ respectively. In particular,
$S'$ and $T'$ are regular operators.
Furthermore, since
$$
N(S')\cap R(T')=\tilde{Y},\hskip1cm N(T')\cap R(S')= \tilde{X},
$$ 
it is clear that
\begin{align*}
&N(S')/(N(S')\cap R(T'))=Y_1,& &R(T')/(N(S')\cap R(T'))=Y_2,\\
&N(T')/(N(T')\cap R(S'))=X_1,& &R(S')/(N(T')\cap R(S'))=X_2.
\end{align*}
\indent Therefore, $(S',T')$ is a regular Fredholm pair and
\vskip.5cm
$$
\hbox{\rm ind}\hbox{ } (S',T')=- \hbox{\rm ind}\hbox{ } (S,T).
$$
\end{proof}
\indent Before ending this section, the perturbation properties of regular
Fredholm pairs are considered. It is clear that the main results of \cite{2}, 
Theorems 3.1 and 3.2, are still true for regular Fredholm pairs. On the other hand, 
thanks to \cite[6.3.4]{5}, also \cite[2.3]{2} remains true for regular Fredholm pairs.\par

\indent In the next section regular Fredholm pairs
will be entirely classified.\par      
\vskip.3cm
\noindent \bf{3. The classification of regular Fredholm pairs}\rm \vskip.3cm

\indent In order to classify regular Fredholm pairs, two sequences
of subspaces are introduced.\par

\newtheorem*{def3.1}{Definition 3.1}
\begin{def3.1} Let $X$ and $Y$ be two Banach spaces and $(S,T)$ belong to 
$RP(X,Y)$. The sequences 
$(R_{S, n})_{n\in {\mathbb N}_0}$ and $(R_{T, n})_{n\in {\mathbb N}_0}$ are defined in
the following way:\par
\indent If $n=0$, then
$$
R_{S,0}=Y,\hskip1cm R_{T,0}=X,
$$
and if $R_{S,n}$ and $R_{T,n}$ are defined, then
$$
R_{S,n+1}=S(R_{T, n}),\hskip1cm R_{T,n+1}=T(R_{S,n}).
$$
\end{def3.1}

\newtheorem*{rem3.2}{Remark 3.2} 
\begin{rem3.2}\rm 
$(R_{S, n})_{n\in {\mathbb N}_0}$ and 
$(R_{T, n})_{n\in {\mathbb N}_0}$ are decreasing sequences of $Y$ and $X$ respectively. 
In fact, 
$$
R_{S,1}=R(S)\subseteq Y=R_{S,0},\hskip1cm R_{T,1}=R(T)\subseteq X=R_{T,0}.
$$
\indent On the other hand, if $R_{S,n}\subseteq R_{S,n-1}$ and $R_{T,n}
\subseteq R_{T,n-1}$,
for a fixed $n\ge 1$, then 
\begin{align*}
&R_{S,n+1}=S(R_{T,n})\subseteq S(R_{T,n-1})=R_{S,n},\\
&R_{T,n+1}=T(R_{S,n})\subseteq T(R_{S,n-1})=R_{T,n}.
\end{align*}
\indent Furthermore, since $R_{S,2}=R(ST)$ and $R_{T,2}=R(TS)$, 
$R_{S,n}$ and $R_{T,n}$ are finite dimensional subspaces of $Y$ and
$X$ respectively, $n\ge 2$.\par
\end{rem3.2}
\indent Next follows a description of the above sequences of subspaces.\par

\newtheorem*{prop3.3}{Proposition 3.3} 
\begin{prop3.3}Let $X$ and $Y$ be two Banach spaces and $(S,T)$ belong to $RP(X,Y)$. 
Then given $n\in {\mathbb N}$
there are subspaces of $X$, $N^n$ and $X_2^n$, and of $Y$, $M^n$ and $Y_2^n$, 
such that: 

i) $R_{S,n}= M^n\oplus Y_2^n$, $R_{T,n}=N^n\oplus X_2^n$;\par
ii) $M^n=R_{S,n}\cap N(T)$, $N^n=R_{T,n}\cap N(S)$;\par
iii) $Y_2^n= R_{S,n}\cap Y_2^k$, $X_2^n=R_{T,n}\cap X_2^k$, $k=1,\ldots n-1$;\par
iv) $(M^n)_{n\in {\mathbb N}}$ and $(N^n)_{n\in{\mathbb N}}$ are decreasing
sequences of subspaces contained in $N(T)$ and $N(S)$ respectively,
moreover, $M^n$ and $N^n$ are finite dimensional subspaces for $n\ge 2$;\par  
v) $(Y_2^n)_{n\in {\mathbb N}}$ and $(X_2^n)_{n\in{\mathbb N}}$ are decreasing
sequences of finite dimensional subspaces contained in $Y_2$ and $X_2$ 
respectively;\par  
vi) $S$ (resp. $T$) induces an isomorphism
$$
X_2^n\xrightarrow{\cong} R_{S,n+1}\hskip1cm (resp.\hbox{  } Y_2^n
\xrightarrow{\cong} R_{T,n+1}).
$$
\end{prop3.3} 

\begin{proof}

\indent When $n=1$ define 
$$
M^1=N(T)\cap R(S),\hbox{  } X_2^1=X_2,\hbox{  }, N^1=N(S)\cap 
R(T),\hbox{  } Y_2^1=Y_2,
$$
where $X_2$ and $Y_2$ are the subspaces considered in Remark 2.2 and 
Remark 2.6.\par
\indent It is clear that these subspaces verify the above assertions. Next suppose that 
the propositon is true for $n\ge 1$. Since $T$ (resp. $S$) induces an isomorphism
\begin{align*} 
&R_{S,n+1}/(N(T)\cap R_{S,n+1})\xrightarrow{\cong} R_{T,n+2},\\  
&(resp. \hbox{  }R_{T,n+1}/(N(S)\cap R_{T,n+1})\xrightarrow{\cong} R_{S,n+2},
\end{align*} 
there are finite dimensional subspaces $V$ and $W$ of $Y$ and $X$ respectively, 
such that
$$
R_{S,n+1}=(N(T)\cap R_{S, n+1})\oplus V, \hskip1cm R_{T,n+1}=(N(S)\cap R_{T,n+1})
\oplus W,
$$
and $T$ (resp. $S$) induces an isomorphism 
$$
V\xrightarrow{\cong} R_{T,n+2} \hskip1cm (\hbox{resp.  } W\xrightarrow{\cong} 
R_{S,n+2}). 
$$
\indent Observe that, according to an argument similar to one used in Remark 2.2,
it is possible to choose $V\subseteq Y_2^n$ and $W\subseteq X_2^n$.
Then, define
$$
M^{n+1}=N(T)\cap R_{S, n+1},\hbox{  }N^{n+1}=N(S)\cap R_{T, n+1},\hbox{  }
Y_2^{n+1}=V,\hbox{ } X_2^{n+1}=W.
$$

\indent Clearly $Y_2^{n+1}\subseteq R_{S, n+1}\cap Y_2^n\subseteq R_{S,n+1}
\cap Y_2^k
\subseteq R_{S,n+1}\cap Y_2$, for $k=1,\ldots ,n$. On the other hand,
if $a\in R_{S, n+1}\cap Y_2$, since $R_{S,n+1}= (N(T)\cap R_{S, n+1})\oplus Y_2^{n+1}$,
then there are $m\in N(T)\cap R_{S, n+1}$ and $y\in Y_2^{n+1}$ such that $a=m+y$. 
However, since $m\in N(T)$ and $a-y\in Y_2$, for $Y_2^{n+1}\subseteq Y_2^n\subseteq 
Y_2$, then $m=0$ and $a=y\in Y_2^{n+1}$. Therefore, $Y_2^{n+1}=R_{S,n+1}\cap 
Y_2^k=R_{S,n+1} \cap Y_2$, for $k=1,\ldots n$. \par
\indent Similarly, $X_2^{n+1}=R_{T,n+1}\cap X_2^k=R_{T,n+1}\cap X_2$, 
for $k=1,\ldots ,n$.\par
\indent The other points of the proposition are clear.\par
\end{proof}

\indent Our next step consists in the description of the relationship between
$R_{S,n}$ and $R_{S,n+1}$, and between $R_{T,n}$ and $R_{T,n+1}$. 
However, to this end, it is necessary to introduce two new sequences of 
subspaces.\par

\newtheorem*{def3.4}{Definition 3.4} 
\begin{def3.4}Let $X$ and $Y$ be two Banach spaces and $(S,T)$ 
belong to $RP(X,Y)$. The sequences of subspaces of $Y$ and 
$X$, $(R_{\tilde{S}, n})_{n\in \mathbb N}$ and $(R_{\tilde{T}, n})_{n\in \mathbb N}$ 
respectively, are defined in the following way.\par
\indent If $n=1$, then
$$
R_{\tilde{S},1}=R(\tilde{S})= S(\tilde{X}),\hskip1cm 
R_{\tilde{T},1}=R(\tilde{T})= T(\tilde{Y}),
$$
where $\tilde{X}$, $\tilde{Y}$, $\tilde{S}$ and $\tilde{T}$ are the spaces and operators 
introduced in Remark 2.6, and if
$n\ge 2$,
$$
R_{\tilde{S},n+1}=S(R_{\tilde{T},n}),\hskip1cm R_{\tilde{T},n+1}=T(R_{\tilde{S},n}). 
$$ 
\indent Observe that $R_{\tilde{S},n}\subseteq R_{S,n}$ and
$R_{\tilde{T},n}\subseteq R_{T,n}$, for $n\in\mathbb N$.
\end{def3.4} 

\indent In the next proposition the sequences introduced in Definition 3.4 are 
characte-
rized.\par

\newtheorem*{prop3.5}{Proposition 3.5} 
\begin{prop3.5}Let $X$ and $Y$ be two Banach spaces and $(S,T)$ belong to 
$RP(X,Y)$. Then there are four
sequences of subspaces, two of $X$, $(\mathbb{N}^n)_{n\in \mathbb N}$ and 
$(\mathbb{X}_2^n)_{n\in \mathbb N}$, and two of $Y$, $(\mathbb{M}^n)_{n\in \mathbb N}$ and 
$(\mathbb{Y}_2^n)_{n\in \mathbb N}$,
such that for $n\in\mathbb N$:\par
i) $R_{\tilde{S}, n}=\mathbb{M}^n\oplus \mathbb{Y}_2^n$, $R_{\tilde{T}, n}=
\mathbb{N}^n\oplus \mathbb{X}_2^n$;\par
ii) $\mathbb{M}^n=R_{\tilde{S},n}\cap N(T)\subseteq M^n$, $\mathbb{N}^n=
R_{\tilde{T},n}\cap N(S)\subseteq N^n$;\par
iii) $\mathbb{Y}_2^n=R_{\tilde{S},n}\cap Y_2^n\subseteq Y_2^n$, $\mathbb{X}_2^n=
R_{\tilde{T},n}\cap X_2^n\subseteq X_2^n$;\par 
iv) $T$ (resp. $S$) induces an isomorphism
$$
\mathbb{Y}_2^n\xrightarrow{\cong}  R_{\tilde{T},n+1},\hskip1cm (resp.\hbox{  } 
\mathbb{X}_2^n  \xrightarrow{\cong}
R_{\tilde{S}, n+1}); 
$$
\indent v) $R_{S,n}=R_{S,n+1}\oplus R_{\tilde{S},n}$, $R_{T,n}=R_{T,n+1}
\oplus R_{\tilde{T},n}$;\par
vi) $M^n=M^{n+1}\oplus \mathbb{M}^n$, $N^n=N^{n+1}\oplus \mathbb{N}^n$;\par
vii) $Y_2^n=Y_2^{n+1}\oplus \mathbb{Y}_2^n$, $X_2^n=X_2^{n+1}\oplus 
\mathbb{X}_2^n$;\par 
viii) when $n=1$ there are subspaces of $\tilde{X}$, $\tilde{X}_N$ and $\tilde{X}_2$,
and of $\tilde{Y}$, $\tilde{Y}_N$ and $\tilde{Y}_2$, such that $\tilde{X}_2$ and
$\tilde{Y}_2$ are finite dimensional,
$$
\tilde{X}=\tilde{X}_N\oplus\tilde{X}_2,\hskip1cm \tilde{Y}=\tilde{Y}_N\oplus\tilde{Y}_2,
$$
and the following operators are isomorphic maps:
$$
\tilde{S}\colon \tilde{X}_N \xrightarrow{\cong} \mathbb{M}^1,\hskip1cm \tilde{S}\colon 
\tilde{X}_2\xrightarrow{\cong} \mathbb{Y}_2^1,
$$
$$
\tilde{T}\colon \tilde{Y}_N \xrightarrow{\cong} \mathbb{N}^1,\hskip1cm \tilde{T}\colon 
\tilde{Y}_2\xrightarrow{\cong}\mathbb{X}_2^1.
$$ 
\end{prop3.5}
 
\begin{proof}

\indent Given $n\in\Bbb N$, consider the isomorphism induced by $T$
$$
R_{\tilde{S},n}/(N(T)\cap R_{\tilde{S}, n})\xrightarrow{\cong} R_{\tilde{T},n+1}.
$$
\indent Since $R_{\tilde{T},n+1}\subseteq R_{T,n+1}\subseteq R_{T,2}=R(TS)$, 
there is $L_n$, a finite dimensional subspace of $R_{\tilde{S},n}$, such that
$R_{\tilde{S},n}=(N(T)\cap R_{\tilde{S},n})\oplus L_n$. 
In addition, $T$ induces an isomorphism
$$
L_n \xrightarrow{\cong} R_{\tilde{T}, n+1}.
$$
\indent Furthermore, since $R_{\tilde{S},n}\subseteq R_{S,n}=M^n\oplus Y_2^n=
(R_{S,n}\cap
N(T))\oplus Y_2^n$, according to an argument similar to the one used in Remark 2.2,
it is possible to choose $L_n\subseteq Y_2^n$. \par
\indent Now well, it is clear that $L_n\subseteq R_{\tilde{S},n}\cap Y_2^n$.
On the other hand, if $a\in R_{\tilde{S},n}\cap Y_2^n$, then there are $m\in N(T)\cap
R_{\tilde{S},n}$ and  $l\in L_n$ such that $a=m+l$. However, since $a-l\in Y_2^n\cap
(N(T)\cap R_{\tilde{S},n})\subseteq Y_2^n\cap (N(T)\cap R_{S,n})=Y_2^n\cap M^n$,
$m=0$ and $a=l\in L_n$. Therefore, $L_n= R_{\tilde{S},n}\cap Y_2^n$.\par 

\indent Next define 
$$
\mathbb{M}^n=R_{\tilde{S},n}\cap N(T),\hskip1cm \mathbb{Y}_2^n=L_n.
$$ 
\indent It is clear that assertions i)-iv) have been proved for $S$.
A similar argument proves the same points for the operator $T$.\par
\indent In order to prove v), an inductive argument wil be used.\par
\indent According to Remark 2.6, 
$$
R_{S,1}=R_{S,2}\oplus R_{\tilde{S},1},\hskip1cm R_{T,1}=R_{T,2}\oplus R_{\tilde{T},1}.
$$
\indent Next suppose that the point v) is true for the operators $S$ and $T$ and for 
a fixed $n\ge 1$. Then, according to Proposition 3.3 i), ii) and vi), and to
Propositon 3.5 i), ii) and iv), which have just been proved,
\begin{align*}
R_{S,n+1}=&S(R_{T,n})=S(R_{T,n+1}\oplus R_{\tilde{T},n})= S((N^{n+1}\oplus X_2^{n+1})
\oplus (\mathbb N^n\oplus \mathbb{X}_2^n))\\
=&S(X_2^{n+1}\oplus \mathbb{X}_2^n)= S(X_2^{n+1})\oplus S(\mathbb{X}_2^n)=R_{S,n+2}
\oplus R_{\tilde{S}, n+1}.
\end{align*}

\indent Similarly, $R_{T,n+1}= R_{T,n+2}\oplus R_{\tilde{T},n+1}$.\par 

\indent Next, according to Proposition 3.3 i) and to Proposition 3.5 i) and v), which 
have just been proved, it is easy to conclude that
$$
M^n\oplus Y_2^n= (M^{n+1}\oplus \mathbb{M}^n)\oplus (Y_2^{n+1}\oplus \mathbb{Y}_2^n),
\hbox{  }
N^n\oplus X_2^n= (N^{n+1}\oplus \mathbb{N}^n)\oplus (X_2^{n+1}\oplus \mathbb{X}_2^n).
$$ 
\indent In particular,
\begin{align*}
&M^n=M^{n+1}\oplus \mathbb{M}^n,&  &N^n=N^{n+1}\oplus \mathbb{N}^n.\\
&Y_2^n=Y_2^{n+1}\oplus \mathbb{Y}_2^n,&  &X_2^n=X_2^{n+1}\oplus 
\mathbb{X}_2^n.
\end{align*}
\indent Finally, consider $n=1$. According to Remark 2.6, 
$$
\tilde{S}\colon \tilde{X}\xrightarrow{\cong}
R_{\tilde{S},1}=\mathbb{M}^1\oplus \mathbb{Y}_2^1.
$$
\indent Therefore, if 
$$
\tilde{X}_N=\tilde{S}^{-1}(\mathbb{M}^1),\hskip1cm \tilde{X}_2=\tilde{S}^{-1}(\mathbb{Y}_2^1),
$$
then $\tilde{X}=\tilde{X}_N\oplus\tilde{X}_2$, 
$$
\tilde{S}\colon \tilde{X}_N \xrightarrow{\cong}\mathbb{M}^1,\hskip1cm \tilde{S}\colon 
\tilde{X}_2\xrightarrow{\cong} \mathbb{Y}_2^1,
$$
and since $\mathbb{Y}_2^1\subseteq Y_2$, $\tilde{X}_2$ is finite dimensional
subspace of $X$.\par
\indent A similar argument proves the case $n=1$ for the operator $T$.\par
\end{proof}
\indent As a result of Propositions 3.3 and 3.5, descriptions of $X$,
$Y$, $S$, and $T$ are obtained.\par  

\newtheorem*{rem3.6}{Remark 3.6}
\begin{rem3.6} \rm Let $X$ and $Y$ be two Banach spaces and $(S,T)$ belong to 
$RP(X,Y)$. If $n\in \mathbb N$, then according to Remark 2.6 and Propositions 3.3 
and 3.5, $X$ and $Y$ may be presented as
\begin{align*}
&X=[X_1\oplus (N^n\oplus\oplus_{i=1}^{n-1}\mathbb{N}^{i})]\oplus [X_2^n\oplus
\oplus_{i=1}^{n-1}
\mathbb{X}^{i}_2]\oplus[\tilde{X}_N\oplus\tilde{X}_2],\\
&Y=[Y_1\oplus (M^n\oplus\oplus_{i=1}^{n-1}\mathbb{M}^{i})]\oplus [Y_2^n
\oplus\oplus_{i=1}^{n-1}
\mathbb{Y}^{i}_2]\oplus[\tilde{Y}_N\oplus\tilde{Y}_2],
\end{align*}
and the maps $S$ and $T$ as
\begin{align*}
&S\mid_{X_1\oplus (N^n\oplus\oplus_{i=1}^{n-1}\mathbb{N}^{i})}\equiv 0,&  
&S\mid_{X_2^n}\colon X_2^n\xrightarrow{\cong}R_{S,n+1}=M^{n+1}\oplus Y_2^{n+1}, \\
&S\mid_{\mathbb{X}_2^{i}}\colon \mathbb{X}_2^{i}\xrightarrow{\cong}R_{\tilde{S},i+1}=
\mathbb{M}^{i+1}\oplus \mathbb{Y}_2^{i+1},& &S\mid_{\tilde{X}_N}\colon \tilde{X}_N
\xrightarrow{\cong} \mathbb{M}^1,\hbox{  }
S\mid_{\tilde{X}_2}\colon \tilde{X}_2\xrightarrow{\cong}\mathbb{Y}_2^1,\\ 
&T\mid_{Y_1\oplus (M^n\oplus\oplus_{i=1}^{n-1}\mathbb{M}^{i})}\equiv 0,& 
&T\mid_{Y_2^n}\colon Y_2^n\xrightarrow{\cong}R_{T,n+1}=N^{n+1}\oplus X_2^{n+1}, \\
&T\mid_{\mathbb{Y}_2^{i}}\colon \mathbb{Y}_2^{i}\xrightarrow{\cong}R_{\tilde{T},i+1}=
\mathbb{N}^{i+1}\oplus \mathbb{X}_2^{i+1},& &T\mid_{\tilde{Y}_N}\colon \tilde{Y}_N
\xrightarrow{\cong} \mathbb{N}^1,\hbox{  }
T\mid_{\tilde{Y}_2}\colon \tilde{Y}_2\xrightarrow{\cong} \mathbb{X}_2^1, 
\end{align*}
where $i=1,\ldots ,n-1$.\par
\end{rem3.6}
\newtheorem*{rem3.7}{Remark 3.7}
\begin{rem3.7} \rm Let $X$ and $Y$ be two Banach spaces and $(S,T)\in RP(X,Y)$. 
Consider the sequences of subspaces of $X$ and $Y$ $(R_{S,n})_{n\in \mathbb N}$ 
and $(R_{T,n})_{n\in \mathbb N}$ respectively. Since 
$R_{S,n}\subseteq R_{S,2}=R(ST)$ for $n\ge 2$, there is $n_0\ge 2$ such that 
$R_{S,n_0}
=R_{S,n_0+1}$. Furthermore, according to this observation, it is easy to prove that
there is $l\in \mathbb N$ such that $R_{S,l}=R_{S,l+k}$ for $k\ge 0$. \par
\indent Similarly, there is $m\in \mathbb N$ such that $R_{T,m}=
R_{T, m+k}$ for $k\ge 0$. \par 
\indent Now well, if $R_{S,l}=R_{S,l+k}$ for $k\ge 0$, then 
$R_{T,l+1}=R_{T,l+1+k}$ for $k\ge 0$. Similarly,
if $R_{T,m}=R_{T,m+k}$ for $k\ge 0$, then $R_{S,m+1}=R_{S,m+1+k}$ for
$k\ge 0$. Therefore, if p and q denote
the first natural numbers such that $R_{S,p}=R_{S,p+k}$ and $R_{T,q}=R_{T,q+k}$
for $k\ge 0$, then $p$, $q\in \mathbb{N}$, and there are the following possibilities:\par
i) $p=q$,\par
ii) if $p<q$, then $q=p+1$,\par
iii) if $q<p$, then $p=q+1$.\par 
\end{rem3.7}
\indent The previous remark leads to a definitions which is central in the classification 
of regular Fredholm pairs.\par
\newtheorem*{def3.8}{Definition 3.8}
\begin{def3.8} Let $X$ and $Y$ be two Banach spaces, $(S,T)\in RP(X,Y)$ and $p$ and 
$q$ as in Remark 3.7.
It will be said that the number of the pair $(S,T)$ is $n$, if $n=min\{ p,q\}$, and it
will be said that the case of the pair $(S,T)$
is  $I-n$ if $p=q$,
 $II-n$ if $p<q$, and $III-n$ if $q<p$.
\end{def3.8} 
\indent Observe that the above construction is symmetric in $X$ and $Y$ and in 
$S$ and $T$. Consequently, in order to study regular Fredholm pairs, interchanging 
$X$ with $Y$ and $S$ with $T$ if necessary, 
it is enough to consider only cases $I-n$ and $II-n$.\par
\indent  In the following theorems the classification of regular Fredholm pairs is 
achieved.
Note that the notations of Remark 3.7 will be used. First of all, regular Fredholm 
pairs whose numbers are equal to $1$ are considered.\par
\newtheorem*{th3.9}{Theorem 3.9} 
\begin{th3.9}Let $X$ and $Y$ be two Banach spaces and $(S,T)$ belong to $RP(X,Y)$. 
Suppose that the number of $(S,T)$ is 1.\par
\indent If the case of $(S,T)$ is 
$I-1$, then the spaces $X$ and $Y$ can be presented as 
$$
X=X_1 \oplus X_2^2,\hskip1cm Y=Y_1 \oplus Y_2^2,
$$ 
and the operators $S$ and $T$ as 
\begin{align*}
&S \mid_{X_1}\equiv 0,& &S\mid_{X_2^2}\colon X_2^2\xrightarrow{\cong}
 Y_2^2,\\
&T \mid_{Y_1}\equiv 0,& &T\mid_{Y_2^2}\colon Y_2^2\xrightarrow{\cong}
 X_2^2.
\end{align*}
\indent If the case of $(S,T)$ is $II-1$, then the spaces $X$ and $Y$ can be
presented as
$$
X=[X_1 \oplus \mathbb{N}^1]\oplus X_2^2,\hskip1cm Y=Y_1 \oplus Y_2^2\oplus 
\tilde{Y}_N,
$$
and the operators $S$ and $T$ as
\begin{align*}
&S \mid_{X_1\oplus \mathbb{N}^1 }\equiv 0,& &S\mid_{X_2^2}\colon X_2^2
\xrightarrow{\cong} Y_2^2,&\\
&T \mid_{Y_1}\equiv 0,& &T\mid_{Y_2^2}\colon Y_2^2\xrightarrow{\cong}
 X_2^2,&
&T\mid_{\tilde{Y}_N}\colon \tilde{Y}_N\xrightarrow{\cong}\mathbb{N}^1.
\end{align*}
\indent If the case of $(S,T)$ is $III-1$, then the spaces $X$ and $Y$ can be
presented as
$$
X=X_1 \oplus X_2^2\oplus\tilde{X}_N,\hskip1cm Y=[Y_1\oplus\mathbb{M}^1] \oplus Y_2^2,
$$
and the operators $S$ and $T$ as
\begin{align*}
& S \mid_{X_1 }\equiv 0,& &S\mid_{X_2^2}\colon X_2^2\xrightarrow{\cong} 
Y_2^2,&
&S\mid_{\tilde{X}_N}\colon \tilde{X}_N\xrightarrow{\cong}\mathbb{M}^1,\\
&T \mid_{Y_1\oplus \mathbb{M}^1}\equiv 0,& &T\mid_{Y_2^2}\colon Y_2^2
\xrightarrow{\cong} X_2^2.
\end{align*}
\indent The spaces involved in the above decompositions are the ones of 
Remark 3.6.\par  
\end{th3.9}
\begin{proof}
\indent Suppose that the case of $(S,T)$ is $I-1$. Since $R(S)=R_{S,1+k}$, 
$k\ge 0$, according
to Proposition 3.5 v), $R_{\tilde{S},k}=0$ for $k\ge 1$.
In particular, $R_{\tilde{S},1}=S(\tilde{X})=0$, and 
since $\tilde{S}\colon \tilde{X}\to R(\tilde{S})$ is an isomorphic map,
according to Proposition 3.5 viii), $\tilde{X}=0$, $\mathbb{M}^1=0$ and $\mathbb{Y}^1_2=0$.
\par
\indent Similarly, since $1=p=q$, $R_{\tilde{T},k}=0$ for $k\ge 1$, $\tilde{Y}=0$, 
$\mathbb{N}^1=0$ and $\mathbb{X}^1_2=0$.
Therefore, according to Remark 3.6, $X$ and $Y$ can be presented as
$$
X=[X_1\oplus N^2] \oplus X_2^2,\hskip1cm Y=[Y_1\oplus M^2 ]\oplus Y_2^2,
$$ 
and $S$ and $T$ as 
\begin{align*}
&S\mid_{X_1\oplus N^2}\equiv 0,&    &S\mid_{X_2^2}\colon X_2^2
\xrightarrow{\cong}M^2 \oplus Y_2^2,\\
&T\mid_{Y_1\oplus M^2}\equiv 0,& &T\mid_{Y_2^2}\colon Y_2^2
\xrightarrow{\cong} N^2 \oplus X_2^2.
\end{align*}
\indent Now well, according to the above presentation, $\dim$ $X_2^2=$ 
$\dim$ $Y_2^2$,
$M^2=0$ and $N^2=0$.\par 
\indent Next suppose that $p=1$ and $q=p+1=2$. According to
Remark 3.7, this is equivalent to the fact that  
$R(S)=R_{S,1+k}$ and $R_{T,2}=R_{T,2+k}$, where $k$ is a positive integer. 
Therefore, according to Proposition 3.5 v), $R_{\tilde{S},1+k}=0$ for $k\ge 0$. 
In particular, and according to
Proposition 3.5 viii), 
$R_{\tilde{S},1}=0$, $\tilde{X}=0$, $\mathbb{M}^1=0$ and $\mathbb{Y}_2^1=0$.
Similarly, and according to Proposition 3.5 v), $R_{\tilde{T},2+k}=0$ for $k\ge 0$. In 
particular, $R_{\tilde{T},2}=0$, and according to Proposition 3.5 vi) and vii),
$\mathbb{N}^2=0$, $\mathbb{X}_2^2=0$, $N^2=N^3$ and $X_2^2=X_2^3$. As a consequence, 
according to Remark 3.6,
the spaces $X$ and $Y$ can be presented as
$$  
X=[X_1\oplus (N^2 \oplus \mathbb{N}^1)]\oplus [X_2^2\oplus \mathbb{X}^1_2],\hskip.5cm
Y=[Y_1\oplus M^2] \oplus Y_2^2\oplus \tilde{Y},
$$ 
and $S$ and $T$ as   
\begin{align*}
&S\mid_{X_1\oplus (N^2\oplus\mathbb{N}^1)}\equiv 0,& &S\mid_{X_2}\colon 
X_2\xrightarrow{\cong}M^2 \oplus Y_2^2,\\
&T\mid_{Y_1\oplus M^2}\equiv 0,& &T\mid_{Y_2^2}\colon Y_2^2
\xrightarrow{\cong}N^2 \oplus X_2^2,& &T\mid_{\tilde{Y}}\colon \tilde{Y}
\xrightarrow{\cong} \mathbb{N}^1\oplus \mathbb{X}_2^1.
\end{align*}
\indent Now well, since $R_{\tilde{S},2}=0$, according to Proposition 3.5 iv),
$\mathbb{X}_2^1=0$. Therefore, $X_2=X_2^2$, $\dim$ $X_2^2$= $\dim$ $Y_2^2$,
$N^2=0$, $M^2=0$, and according to Proposition 3.5 viii), $\tilde{Y}_2=0$.
\par
\indent Interchanging $X$ with $Y$ and $S$ with $T$, the proof of the case $III-1$
can be carried out with an argument similar to the one of the case $II-1$.\par 
\end{proof}
\indent Observe that in the case $I-1$, $X$ and $Y$ are always finite dimensional 
Banach spaces.\par
\indent Next, regular Fredholm pairs whose numbers are greater or equal to $2$
are classified.\par  
\newtheorem*{th3.10}{Theorem 3.10} 
\begin{th3.10} Let $X$ and $Y$ be two Banach spaces and $(S,T)$ belong to $RP(X,Y)$. Suppose that the case of $(S,T)$ is 
$I-n$. Then, If $n=p=q\ge2$ is the number of the pair $(S,T)$, the spaces $X$ and $Y$
can be presented as 
\begin{align*}
&X=[X_1 \oplus\oplus_{i=1}^{p-1}\mathbb{N}^{i}]\oplus [X_2^{p-1}\oplus\oplus_{i=1}^{p-2}
\mathbb{X}^{i}_2]\oplus[\tilde{X}_N\oplus\tilde{X}_2],\\
&Y=[Y_1 \oplus\oplus_{i=1}^{p-1}\mathbb{M}^{i}]\oplus [Y_2^{p-1}
\oplus\oplus_{i=1}^{p-2}\mathbb{Y}^{i}_2]\oplus[\tilde{Y}_N\oplus\tilde{Y}_2],
\end{align*}
and the operators $S$ and $T$ as
\begin{align*}
&S\mid_{X_1\oplus \oplus_{i=1}^{p-1}\mathbb{N}^{i}}\equiv 0,& &S\mid_{X_2^{p-1}}
\colon X_2^{p-1}\xrightarrow{\cong} Y_2^{p-1},& &S\mid_{\mathbb{X}_2^{i}}\colon 
\mathbb{X}_2^{i}\xrightarrow{\cong}\mathbb{M}^{i+1}\oplus \mathbb{Y}_2^{i+1},\\
&S\mid_{\mathbb{X}_2^{p-2}}\colon \mathbb{X}_2^{p-2}\xrightarrow{\cong}
\mathbb{M}^{p-1},&
&S\mid_{\tilde{X}_N}\colon \tilde{X}_N\xrightarrow{\cong} \mathbb{M}^1,&
&S\mid_{\tilde{X}_2}\colon \tilde{X}_2\xrightarrow{\cong} \mathbb{Y}_2^1,\\ 
&T\mid_{Y_1\oplus\oplus_{i=1}^{p-1}\mathbb{M}^{i}}\equiv 0,& &T\mid_{Y_2^{p-1}}
\colon Y_2^{p-1}\xrightarrow{\cong} X_2^{p-1},&
&T\mid_{\mathbb{Y}_2^{i}}\colon \mathbb{Y}_2^{i}\xrightarrow{\cong}\mathbb{N}^{i+1}
\oplus \mathbb{X}_2^{i+1},\\ 
&T\mid_{\mathbb{Y}_2^{p-2}}\colon \mathbb{Y}_2^{p-2}\xrightarrow{\cong}
\mathbb{N}^{p-1},&
&T\mid_{\tilde{Y}_N}\colon \tilde{Y}_N\xrightarrow{\cong} \mathbb{N}^1,&
&T\mid_{\tilde{Y}_2}\colon \tilde{Y}_2\xrightarrow{\cong} \mathbb{X}_2^1, 
\end{align*}
where $i=1,\ldots ,p-3$, and the spaces involved in the above decomposition 
are the ones of
Remark 3.6.\par
\indent In addition, if $n=2$, then $\mathbb{X}_2^1$, $\mathbb{Y}_2^1$, 
$\mathbb{Y}_2$ and $\mathbb{X}_2$ are null spaces.\par
\end{th3.10}
\begin{proof} 
\indent Let $p=q\ge 2$ be the number of the pair $(S,T)$. Since $R_{S,p}=R_{S,p+k}$
and $R_{S,p+k}=R_{S, p+k+1}\oplus R_{\tilde{S}, p+k}$ for $k\ge 0$, then 
$R_{\tilde{S},p+k}=0$
for $k\ge 0$, that is $\mathbb{M}^{p+k}=0$ and $\mathbb{Y}_2^{p+k}=0$ for $k\ge 0$. 
Furthermore,
since according to Proposition 3.5 iv) and vii),
$$
S\colon\Bbb{X}_2^{p-1+k}\xrightarrow{\cong} R_{\tilde{S}, p+k},\hskip1cm 
X_2^{p-1+k}=X_2^{p+k}\oplus\Bbb{X}_2^{p-1+k},
$$ 
then $\mathbb{X}_2^{p-1+k}=0$ and $X_2^{p-1}=X_2^{p-1+k}$ for $k\ge 0$.
In particular, according to Proposition 3.5 vii),   
$$
X_2=X_2^{p-1}\oplus\oplus_{i=1}^{p-2}\mathbb{X}^{i}_2.
$$
\indent On the other hand, since $p=q$, similar properties can be obtained for 
$T$ and $X$. 
Therefore, $\mathbb{N}^{p+k}=0$, $\mathbb{X}_2^{p+k}=0$, 
$\mathbb{Y}_2^{p-1+k}=0$, $Y_2^{p-1}=Y_2^{p-1+k}$ for $k\ge 0$, and 
$$
Y_2=Y_2^{p-1}\oplus\oplus_{i=1}^{p-2} \mathbb{Y}^{i}_2.
$$
\indent In addition, since according to Proposition 3.3 i) and vi),
\begin{align*}
&S\colon X_2^{p-1}\xrightarrow{\cong} R_{S,p}=M^p\oplus Y_2^p= M^p\oplus Y_2^{p-1},
\\
&T\colon Y_2^{p-1}\xrightarrow{\cong}R_{T,p}=N^p\oplus X_2^p= N^p\oplus X_2^{p-1},
\end{align*}
it is clear that $M^p=0$, $N^p=0$ and
$$
S\colon X_2^{p-1}\xrightarrow{\cong} Y_2^{p-1},\hskip1cm T\colon Y_2^{p-1}
\xrightarrow{\cong} X_2^{p-1}.
$$
\indent Finally, if $p=q=2$, then $R_{\tilde{S},2}=0$ and $R_{\tilde{T},2}=0$.
Consequently, according to Proposition 3.5 i) and viii), $\mathbb{X}_2^1$, 
$\mathbb{Y}_2^1$, $\tilde{X}_2$ and $\tilde{Y}_2$ are null spaces.\par
\end{proof}

\newtheorem*{th3.11}{Theorem 3.11} 
\begin{th3.11}
Let $X$ and $Y$ be two Banach spaces and $(S,T)$ belong to $RP(X,Y)$. Suppose 
that the case of $(S,T)$ is 
$II-n$. Then, if $n=p=q-1\ge2$ is the number of the pair $(S,T)$, the spaces $X$ and 
$Y$ can be presented as
\begin{align*}
&X=[X_1 \oplus\oplus_{i=1}^p\Bbb{N}^{i}]\oplus [X_2^{p-1}\oplus\oplus_{i=1}^{p-2}
\mathbb{X}^{i}_2]\oplus[\tilde{X}_N\oplus\tilde{X}_2],\\
&Y=[Y_1 \oplus\oplus_{j=1}^{p-1}\mathbb{M}^{j}]\oplus [Y_2^p\oplus\oplus_{j=1}^{p-1}
\mathbb{Y}^{j}_2]\oplus[\tilde{Y}_N\oplus\tilde{Y}_2],
\end{align*}
and the operators $S$ and $T$ as
\begin{align*}
&S\mid_{X_1\oplus \oplus_{i=1}^p\Bbb{N}^{i}}\equiv 0,& &S\mid_{X_2^{p-1}}
\colon X_2^{p-1}\xrightarrow{\cong} Y_2^p,&
&S\mid_{\mathbb{X}_2^{i}}\colon \mathbb{X}_2^{i}\xrightarrow{\cong}\mathbb{M}^{i+1}
\oplus \mathbb{Y}_2^{i+1},\\ 
&S\mid_{\tilde{X}_N}\colon \tilde{X}_N\xrightarrow{\cong} 
\mathbb{M}^1,& &S\mid_{\tilde{X}_2}\colon \tilde{X}_2\xrightarrow{\cong} \mathbb{Y}_2^1,\\ 
&T\mid_{Y_1\oplus\oplus_{j=1}^{p-1}\mathbb{M}^{j}}\equiv 0,& &T\mid_{Y_2^p}
\colon Y_2^p\xrightarrow{\cong} X_2^{p-1},& &T\mid_{\mathbb{Y}_2^{j}}\colon 
\mathbb{Y}_2^{j}\xrightarrow{\cong}\mathbb{N}^{j+1}\oplus 
\mathbb{X}_2^{j+1},\\ 
&T\mid_{\mathbb{Y}_2^{k}}\colon \mathbb{Y}_2^{k}\xrightarrow{\cong}\mathbb{N}^{k+1},&
&T\mid_{\tilde{Y}_N}\colon \tilde{Y}_N \xrightarrow{\cong} \mathbb{N}^1,&
&T\mid_{\tilde{Y}_2}\colon \tilde{Y}_2 \xrightarrow{\cong}\mathbb{X}_2^1,\\ 
\end{align*}
where $i=1,\ldots ,p-2$, $j=1,\ldots ,p-3$, $k=p-2$, $p-1$, and the spaces involved in 
the above decomposition are the ones of Remark 3.6.\par
\indent In addition, if $p=2$, then $\mathbb{X}_2^1$ and $\tilde{Y}_2$ are null spaces.\par
\end{th3.11}
\begin{proof}
\indent Let $p=q-1\ge 2$ be the number of the pair $(S,T)$. Consequently, $q=p+1$, and since
$R_{S,p}=R_{S,p+k}$ and $R_{T,p+1}=R_{T,p+1+k}$ for $k\ge 0$, as in 
Theorem 3.10, $R_{\tilde{S},p+k}=0$ and $R_{\tilde{T},p+1+k}=0$ for $k\ge 0$.
In particular, according to Proposition 3.5 i),
$$
\mathbb{M}^{p+k}=0, \hskip1cm\Bbb{N}^{p+1+k}=0,
$$
for $k\ge 0$.\par
\indent On the other hand, according to Proposition 3.5 iv) and vii), 
$$S\colon\Bbb{X}_2^{p-1+k}\xrightarrow{\cong}R_{\tilde{S},p+k},\hskip1cm
T\colon\Bbb{Y}_2^{p+k}\xrightarrow{\cong} R_{\tilde{T},p+k+1},
$$ 
and  
$$X_2^{p-1+k}=X_2^{p+k}\oplus\Bbb{X}_2^{p-1+k},\hskip1cm Y_2^{p+k}=
Y_2^{p+k+1}\oplus\Bbb{Y}_2^{p+k}.
$$ 
Consequently,
$$
\mathbb{X}_2^{p-1+k}=0,\hbox{  }\mathbb{Y}_2^{p+k}=0,\hbox{  } X_2^{p-1}=X_2^{p-1+k},
\hbox{  } Y_2^p=Y_2^{p+k},
$$ 
for $k\ge 0$. Therefore,
$$
X_2=X_2^{p-1}\oplus\oplus_{i=1}^{p-2}\mathbb{X}^{i}_2,\hskip1cm Y_2=Y_2^p\oplus
\oplus_{j=1}^{p-1}\mathbb{Y}^{j}_2.
$$
\indent Now well, since $X_2^{p-1}=X_2^{p+1}$ and since
$$
S\colon X_2^{p-1}\xrightarrow{\cong}R_{S,p}=M^p\oplus Y_2^p,\hbox{  }
T\colon Y_2^p \xrightarrow{\cong}R_{T,p+1}=N^{p+1}\oplus X_2^{p+1},
$$
it is clear that $M^p=0$, $N^{p+1}=0$, and
$$
S\colon X_2^{p-1}\xrightarrow{\cong} Y_2^p,\hskip1cm T\colon Y_2^p
\xrightarrow{\cong} X_2^{p-1}.
$$
\indent Finally, if $p=2$, then $R_{\tilde{S},2}=0$. Consequently, according to
Proposition 3.5 i) and viii), $\mathbb{X}_2^1$ and $\tilde{Y}_2$ are null spaces.\par
\end{proof}

\newtheorem*{th3.12}{Theorem 3.12}
\begin{th3.12} Let $X$ and $Y$ be two Banach spaces and $(S,T)$ belong to $RP(X,Y)$. Suppose that the case of $(S,T)$ is 
$III-n$. Then, if $q=p-1\ge2$ is the number of the pair $(S,T)$, the spaces $X$ and $Y$
can be presented as
\begin{align*}
&X=[X_1 \oplus\oplus_{i=1}^{q-1}\mathbb{N}^{i}]\oplus [X_2^q\oplus\oplus_{i=1}^{q-1}
\mathbb{X}^{i}_2]\oplus[\tilde{X}_N\oplus\tilde{X}_2],\\
&Y=[Y_1 \oplus\oplus_{j=1}^q\Bbb{M}^{j}]\oplus [Y_2^{q-1}\oplus\oplus_{j=1}^{q-2}
\mathbb{Y}^{j}_2]\oplus[\tilde{Y}_N\oplus\tilde{Y}_2],
\end{align*}
and the operators $S$ and $T$ as 
\begin{align*}
&S\mid_{X_1\oplus \oplus_{i=1}^{q-1}\mathbb{N}^{i}}\equiv 0,& &S\mid_{X_2^q}\colon 
X_2^q\xrightarrow{\cong} Y_2^{q-1},& 
&S\mid_{\mathbb{X}_2^{i}}\colon \mathbb{X}_2^{i}\xrightarrow{\cong}\mathbb{M}^{i+1}
\oplus \mathbb{Y}_2^{i+1},\\
&S\mid_{\mathbb{X}_2^{k}}\colon \mathbb{X}_2^{k}\xrightarrow{\cong}\mathbb{M}^{k+1},&
&S\mid_{\tilde{X}_N}\colon \tilde{X}_N\xrightarrow{\cong} \mathbb{M}^1,&
&S\mid_{\tilde{X}_2}\colon \tilde{X}_2\xrightarrow{\cong} \mathbb{Y}_2^1,\\ 
&T\mid_{Y_1\oplus\oplus_{j=1}^q\Bbb{M}^{j}}\equiv 0,& &T\mid_{Y_2^{q-1}}
\colon Y_2^{q-1}\xrightarrow{\cong} X_2^q, &
&T\mid_{\mathbb{Y}_2^{j}}\colon \mathbb{Y}_2^{j}\xrightarrow{\cong}\mathbb{N}^{j+1}\oplus 
\mathbb{X}_2^{j+1},\\ 
&T\mid_{\tilde{Y}_N}\colon \tilde{Y}_N\xrightarrow{\cong} 
\mathbb{N}^1,&
&T\mid_{\tilde{Y}_2}\colon \tilde{Y}_2\xrightarrow{\cong} \mathbb{X}_2^1, 
\end{align*}
where $i=1,\ldots ,q-3$, $k=q-2$, $q-1$, $j=1,\ldots ,q-2$, and the spaces involved 
in the above decomposition are the ones of Remark 3.6.\par
\indent In addition, if $q=2$, then $\mathbb{Y}_2^1$ and $\tilde{X}_2$ are null spaces.\par
\end{th3.12}
\begin{proof}
\indent Interchanging $X$ with $Y$ and $S$ with $T$, the proof can be carried 
out with an 
argument similar to the one 
of Theorem 3.11.\par
\end{proof}
\newtheorem*{rem3.13}{Remark 3.13} 
\begin{rem3.13}\rm Let $X$ and $Y$ be Banach spaces and $(S,T)\in RP(X,Y)$. 
Suppose that the case of $(S,T)$ is $I-2$. According to Theorem 3.10, $X$ and $Y$ 
may be described as
$$
X=(X_1\oplus \mathbb{N}^1)\oplus X_2^1\oplus \tilde{X}_N,\hskip.5cm 
Y=(Y_1\oplus \mathbb{M}^1)\oplus Y_2^1\oplus \tilde{Y}_N,
$$
and the the operatos $S$ and $T$ can be presented as
\begin{align*}
&S\mid_{X_1\oplus \mathbb{N}^1}\equiv 0,& &S\colon X_2^1\xrightarrow{\cong} 
Y_2^1,& &S\colon \tilde{X}_N\xrightarrow{\cong} \mathbb{M}^1,\\
&T\mid_{Y_1\oplus \mathbb{M}^1}\equiv 0,& &T\colon Y_2^1\xrightarrow{\cong} 
X_2^1,&
&T\colon \tilde{Y}_N\xrightarrow{\cong} \mathbb{N}^1.
\end{align*}
\indent Moreover, since $\dim X_2=\dim X_2^1=\dim Y_2^1=\dim Y_2$, 
ind $(S,T)= \dim X_1 -\dim Y_1$. \par

\indent Now well, if $X'$ and $Y'$ are Banach spaces and if $(S',T')\in RP(X',Y')$ 
is a regular Fredholm symmetrical pair, then it is not difficult to prove
that 
$$
X'=X_1'\oplus N\oplus\tilde{X},\hskip1cm Y'=Y_1'\oplus M\oplus\tilde{Y},
$$
where $X_1'$ and $Y_1'$ are finite dimensional subspaces. Furthermore,
the operators $S'$ and $T'$ are such that
$$
S'\mid_{X_1'\oplus N}\equiv 0,\hbox{ } S'\colon \tilde{X}\xrightarrow{\cong} M,\hbox{  }
T'\mid_{Y_1'\oplus M}\equiv 0,\hbox{ } T'\colon \tilde{Y}\xrightarrow{\cong} N,
$$
which implies that ind $(S',T')= \dim X_1-\dim Y_1$. Therefore, a regular Fredholm 
symmetrical pair is nothing but a very particular
type of regular Fredholm pair, that is a pair whose case is $I-2$ and such that
$X_2=0 $ and $Y_2=0$.\par 
\end{rem3.13}
\newtheorem*{rem3.14}{Remark 3.14}
\begin{rem3.14} \rm Observe that if $X$ and $Y$ are Banach
spaces and $(S,T)$ belongs to $RP(X,Y)$, then, according to Theorems 3.9 - 3.12, 
$X$ is a finite dimensional Banach space if and only if $Y$ is.\par
\indent On the other hand, if $X$, $Y$, $S$ and $T$ are constructed as in Theorems
3.9 - 3.12, then $(S,T)\in RP(X,Y)$ and the number and case of $(S,T)$ are the ones 
considered in the corresponding theorem. Therefore, thanks to 
Theorems 3.9 - 3.12, regular Fredholm are entirely classified.\par   
\end{rem3.14}
\vskip.3cm
\noindent \bf{4. The index of regular Fredholm pairs and Weyl Pairs}\rm 
\vskip.3cm
 
\indent In this section the index of a regular Fredholm pair is studied. It is proved that
the index provides a fundamental tool in the description of the spaces and maps
of such a pair. Furthermore, Weyl pairs, that is Fredholm pairs whose index is null,
are considered.\par

\newtheorem*{th4.1}{Theorem 4.1}
\begin{th4.1} Let $X$ and $Y$ be two Banach spaces and $(S,T)$ belong to
$RP(X,Y)$. Suppose that the number of the pair $(S,T)$ is greater or equal $2$. 
Consider $n=2$ and the corresponding decomposition of the spaces $X$ and $Y$
given in Remark 3.6, that is
\begin{align*}
&X=[X_1\oplus (N^2\oplus \mathbb{N}^1)]\oplus [X_2^2\oplus \mathbb{X}^1_2]\oplus
 [\tilde{X}_N\oplus
\tilde{X}_2],\\
&Y=[Y_1\oplus (M^2\oplus \mathbb{M}^1)]\oplus [Y_2^2\oplus \mathbb{Y}^1_2]\oplus 
[\tilde{Y}_N\oplus
\tilde{Y}_2].
\end{align*}
Then
\begin{align*}
\hbox{\rm{ind}}\hbox{  } (S,T)&=\dim\hbox{  } (X_1\oplus N^2) -\dim\hbox{  } (Y_1
\oplus \tilde{Y}_2)\\
                     &=\dim\hbox{  } (X_1\oplus \tilde{X}_2) -\dim\hbox{  } (Y_1\oplus M^2).
\end{align*}
\indent In addition, if the number of the pair $(S,T)$ is $1$, then 
$$
\hbox{\rm{ind}}\hbox{  }(S,T)=\dim\hbox{  }X_1-\dim\hbox{  }Y_1.
$$
\end{th4.1}
\begin{proof}
\indent According to Definition 2.1 and Remark 2.2, the index of the pair $(S,T)$ is 
the number
$$
\hbox{\rm{ind}}\hbox{  } (S,T) =\dim X_1 -\dim X_2 -\dim Y_1+\dim Y_2.
$$
\indent Now well, according to Proposition 3.3 vi), $T\colon Y_2\xrightarrow{\cong}
 R_{T,2}=N^2\oplus X_2^2$. Moreover, according to Proposition 3.5 viii), 
$T\colon \tilde{Y}_2\xrightarrow{\cong} \mathbb{X}^1_2$. Consequently
$$
\dim Y_2-\dim X_2=\dim N^2 -\dim \mathbb{X}^1_2=\dim N^2-\dim \tilde{Y}_2,
$$
and 
$$
\hbox{\rm ind}\hbox{  } (S,T)= \dim (X_1\oplus N^2) -\dim (Y_1\oplus \tilde{Y}_2).
$$
\indent Since ind $(T,S)=-$ ind $(S,T)$, a similar argument proves the second
equality.\par
\indent The last assertion is a consequence of Theorem 3.9.\par
\end{proof}

\newtheorem*{rem4.2}{Remark 4.2}
\begin{rem4.2} \rm Let $X$ and $Y$ be two Banach spaces and $(S,T)$ belong to 
$RP(X,Y)$. Suppose that the number of the pair $(S,T)$ is greater or equal 2. Consider 
again, as in Theorem 4.1, $n=2$ and the corresponding description of the spaces 
$X$ and $Y$ given in Remark 3.6. \par
\indent Now well, according to Proposition 3.5 vii) and viii), $Y_2=Y_2^2\oplus 
\mathbb{Y}^1_2$, and
$$
S\colon \tilde{X}_2\xrightarrow{\cong} \mathbb{Y}^1_2,\hskip.5cm S\colon \tilde{X}_N
\xrightarrow{\cong}\mathbb{M}^1, 
\hskip.5cmT\colon \tilde{Y}_N\xrightarrow{\cong} \mathbb{N}^1.
$$
\indent In addition, according to Proposition 3.3 vi), 
$$S\colon X_2\xrightarrow{\cong} R_{S,2}=M^2\oplus Y_2^2.
$$ 
\indent Consequently, the subspaces
of $X$ and $Y$ that in the above presentation are not related by isomorphic maps 
are $X_1\oplus N^2$ and $Y_1\oplus \tilde{Y}_2$ respectively.\par
\indent Similarly, interchanging $X$ with $Y$ and $S$ with $T$, the subspaces of 
$X$ and $Y$ that in the above presentation are not related  by isomorphic maps 
are $Y_1\oplus M^2$ and $X_1\oplus \tilde{X}_2$ respectively. \par
\indent On the other hand, if the number of the pair is $1$, according to
Theorem 3.9, the subspaces of $X$ and $Y$ that are not related by isomorphic
maps are 
$X_1$ and $Y_1$.\par
\indent As a result, the index has a fundamental role in the description of regular 
Fredholm pairs. In fact, the index is a measure of the subspaces of $X$ and $Y$ 
that in the above decomposition are not related by isomorphisms.\par
\end{rem4.2}

\newtheorem*{rem4.3}{Remark 4.3}
\begin{rem4.3} \rm Let $X$ and $Y$ be two Banach spaces and $(S,T)$ belong to 
$RP(X,Y)$. Suppose that the number of the pair $(S,T)$ is greater or
equal 2. Consider, as in Theorem 4.1, $n=2$ and the corresponding decompositon 
of $X$ and $Y$ given in Remark 3.6. Next suppose that ind $(S,T)=0$. According to 
Theorem 4.1, this is equivalent to the fact that
$X_1\oplus N^2$ is isomorphic to $Y_1\oplus \tilde{Y}_2$
and $X_1\oplus \tilde{X}_2$ to $Y_1\oplus M^2$. However,
since according to Proposition 3.5 viii), $\mathbb{N}^1$ is isomorphic to 
$\tilde{Y}_N$ and 
$\mathbb{M}^1$ to $\tilde{X}_N$, then
$N(S)$ is isomorphic $Y/R(S)$ and $N(T)$ to 
$X/R(T)$. Consequently, according to \cite[3.8.6]{5}, $S$ and $T$ are 
\it{decomposably regular or relatively Weyl operators}\rm, that is $S$ and $T$ are 
regular maps which have isomorphic pseudoinverses $S'\in L(Y,X)$ and
$T'\in L(X,Y)$ respectively, see \cite[3.8.5]{5}. Similarly, if the
number of the pair $(S,T)$ is 1, and if ind $(S,T)=0$, then, according to Theorems 3.9
and 4.1, $S$ and $T$ are decomposably regular operators. As an analogy to Weyl
operators, Weyl pairs are introduced, see \cite[6.5]{5}. \par
\end{rem4.3}
\newtheorem*{def4.4}{Definition 4.4} 
\begin{def4.4} Let $X$ and $Y$ be two Banach spaces and $(S,T)$ belong
to $P(X,Y)$. The pair $(S,T)$ is said a Weyl pair, if \rm{ind}\it\hbox{ }$(S,T)=0$. 
The set of all Weyl pairs is denoted by $W(X,Y)$. In addition, if $(S,T)$ belongs to 
$W(X,Y)\cap RP(X,Y)$, then $(S,T)$ is called a regular Weyl pair. The set of all 
regular Weyl pairs is denoted by $RW(X,Y)$.
\end{def4.4}
\newtheorem*{rem4.5}{Remark 4.5}
\begin{rem4.5}\rm Let $X$ and $Y$ be two Banach spaces
and $S\in L(X,Y)$. According to Remark 2.2, it is clear that if $S$ is a Weyl operator,
then $(S,0)$ belongs to $RW(X,Y)$.\par
\indent On the other hand, it is well known that in order for $S\in L(X,Y)$ to be a Weyl
operator it is necessary and sufficient that $S$ is Fredholm and decomposably regular,
see for example \cite[6.5.2]{5}. However, as the following example shows, there are regular
Fredholm pairs whose operators are decomposably regular and whose index is not 
null.\par

\indent Let $I$ and $J$ arbitrary infinite disjoint sets and consider the Hilbert spaces
$N=l^2(I)$ and $M=l^2(J)$. Let $I_1$ be a finite non void set such that $I_1\cap I=
\emptyset=I_1\cap J$, and consider the Hilbert space
$X_1=l^2(I_1)$. Next define the Hilbert spaces
$$
X=X_1\oplus N\oplus M,\hskip1cm Y=M\oplus N,
$$
and the operators $S\in L(X,Y)$ and $T\in L(Y,X)$
$$
S\mid_{X_1\oplus N}\equiv 0,\hbox{  } S=I_M\colon M\to M,\hbox{ } T\mid_{M} \equiv 0,
\hbox{  } T=I_N\colon N\to N,
$$
where $I_M$ and $I_N$ denotes the identity maps of $M$ and $N$ respectively.\par   
\indent It is clear that $S$ and $T$ are regular operators. Moreover, $(S,T)$ belongs
to $RP(X,Y)$,
actually, $(S,T)$ is a regular Fredholm symmetrical pair, and ind $(S,T)=\dim X_1$, 
which is non null, for $I_1$ is a non void set.\par
\indent Now well, since $I_1$ is a finite set and $I$ is an infinite set, 
$X_1\oplus N$ is isomorphic to $N$. Consequently,
$N(S)$ is isomorphic to $Y/R(S)$, that is $S$ is a decomposably regular operator,
see \cite[3.8.6]{5}. Similarly, $X_1\oplus M$ is isomorphic to
$M$. Therefore, $N(T)$ is isomorphic to $X/R(T)$, that is $T$ is a decomposably
regular operator, see \cite[3.8.6]{5}.
\end{rem4.5}
\vskip.3cm
\noindent \bf{5. Characterizations of Regular Fredholm Pairs}\rm 
\vskip.3cm

\indent In this section three characterizations of regular Fredholm pairs are proved.
In the first one such objects are characterized in terms of regular Fredholm 
symmetrical pairs.  This characterization plays a central role in the proof of the 
second one, where regular Fredholm pairs are characterized in terms of exact 
chains of multiplication operators. Finally, in the third one, 
the objects under consideration are characterized in terms of invertible Banach 
space operators.  \par

\indent In order to prove the first characterization some preparation is needed.\par

\newtheorem*{rem5.1}{Remark 5.1}
\begin{rem5.1} \rm Consider $X$ and $Y$ two Banach spaces, and
$S\in L(X,Y)$ and $T\in L(Y,X)$ two operators such that $R(ST)$ and $R(TS)$ 
are finite dimensional subspaces of $Y$ and $X$ respectively. Then, it is possible 
to define the Banach spaces
$\mathcal{X}=X/R(TS)$ and $\mathcal{Y}=Y/R(ST)$, and the linear and continuous 
maps $\overline{S}\in L(\mathcal{X},\mathcal{Y})$
and $\overline{T}\in L(\mathcal{Y} ,\mathcal{X})$, the factorizations of $S$ and $T$ 
through the respective invariant subspaces. It is clear that $\overline{S}\circ
\overline{T}=0$ and $\overline{T}\circ\overline{S}=0$. Furthermore, according to 
\cite[2.1]{2},
$$
N(\overline{S})/R(\overline{T})\cong N(S)/(N(S)\cap R(T)), \hbox{  }
N(\overline{T})/R(\overline{S})\cong N(T)/(N(T)\cap R(S)).
$$
Therefore, the pair $(S,T)$ belongs to $P(X,Y)$ if and only if $(\overline{S},\overline{T})$
is a Fredholm symmetrical pair, see \cite[2.1]{2}.
\end{rem5.1}
\newtheorem*{th5.2}{Theorem 5.2} 
\begin{th5.2}Let $X$ and $Y$ be two Banach spaces and $(S,T)$ belong
to P(X,Y). Then, with the notations of Remark 5.1, $(S,T)$ belongs to $RP(X,Y)$ 
if and only if $(\overline{S},\overline{T})$
is a regular Fredholm symmetrical pair.
\end{th5.2}
\begin{proof}
\indent First of all, note that if $(S,T)\in RP(X,Y)$, then according to Remark 5.1,
in order to prove that $(\overline{S},\overline{T})$
is a regular Fredholm symmetrical pair, it is enough to prove that $\overline{S}$ and 
$\overline{T}$ are regular operators.\par
 \indent Consider $n=2$ and the corresponding decomposition of $X$ and $Y$ of
Remark 3.6, that is
$$
X=(X_1\oplus N^2\oplus\Bbb{N}^1)\oplus (X_2^2\oplus \mathbb{X}_2^1)\oplus\tilde{X},
\hbox{  } Y=(Y_1\oplus M^2\oplus\Bbb{M}^1)\oplus (Y_2^2\oplus \mathbb{Y}_2^1)\oplus
\tilde{Y},
$$
and recall that 
$$
R(ST)= R_{S,2}=M^2\oplus Y_2^2,\hskip1cm R(TS)= R_{T,2}=N^2\oplus X_2^2.
$$
Therefore, $X/R(TS)$ and $Y/R(ST)$ can be identified with
$$
\mathcal{X}= X_1\oplus \mathbb{N}^1\oplus \mathbb{X}_2^1\oplus\tilde{X},\hskip1cm
\mathcal{Y}= Y_1\oplus \mathbb{M}^1\oplus \mathbb{Y}_2^1\oplus\tilde{Y}.
$$
\indent Moreover, since 
\begin{align*}
&R(\overline{S})= R(S)/R(ST),\hbox{   } R(\overline{T})= R(T)/R(TS),\\
&N(\overline{S})= S^{-1}(R(ST))/R(TS)= (N(S)+ R(T))/R(TS),\\
&N(\overline{T})= T^{-1}(R(TS))/R(ST)= (N(T)+ R(S))/R(ST),
\end{align*}
these spaces can be identified with
\begin{align*}
&R(\overline{S})=\mathbb{M}^1\oplus\Bbb{Y}_2^1,&
&N(\overline{S})=X_1\oplus \mathbb{N}^1\oplus \mathbb{X}_2^1\\
&R(\overline{T})=\mathbb{N}^1\oplus\Bbb{X}_2^1,&
&N(\overline{T})=Y_1\oplus\Bbb{M}^1\oplus \mathbb{Y}_2^1.
\end{align*}
\indent Therefore, according to \cite[3.8.2]{5}, $\overline{S}$ and $\overline{T}$ are regular
operators.\par
\indent On the other hand, suppose that $\overline{S}$ is a regular operator. 
Then, there is $V$, a closed linear subspace of $\mathcal{X}$, such that
$$
N(\overline{S})\oplus V=\mathcal{X}.
$$ 
\indent Let $\pi\colon X\to \mathcal{X}$ be the 
canonical projection and $V_1=\pi^{-1}(V)\cap R(TS)$. Since $V_1$
is a finite dimensional subspace of the Banach space $\pi^{-1}(V)$, there
is a closed linear subspace $W_1\subseteq \pi^{-1}(V)$ such that
$$
V_1\oplus W_1=\pi^{-1}(V).
$$
\indent Now well, since $\pi$ is a surjective map and since $\pi (V_1)=0$,
$$\pi (W_1)=\pi(\pi^{-1}(V))=V.
$$
\indent Furthermore, according to \cite[2.1]{2},
$$
\pi (N(S)+R(T)+W_1)= \pi (N(S) +R(T)) + \pi (W_1)=N(\overline{S}) +V=\mathcal{X}.
$$ 
\indent Consequently, $N(S) + R(T) + W_1+R(TS)=X$. However, since $R(TS)\subseteq
R(T)$, $(N(S)+R(T))+W_1=X$.\par
\indent Next consider $L=(N(S)+R(T))\cap W_1$. Since $\pi (L)\subseteq 
N(\overline{S})\cap V=0$,
$L\subseteq R(TS)$. Consequently, $L\subseteq W_1\cap R(TS) =0$. Therefore,
$$
(N(S)+R(T)) \oplus W_1=X.
$$
\indent Similarly, if $\overline{T}$ is a regular operator, then there is $W_2$, a 
closed subspace of Y, such that 
$$
(N(T)+R(S))\oplus W_2=Y.
$$
\indent Finally, since $(S,T)\in P(X,Y)$, according to Proposition 2.4, $(S,T)$ is a 
regular Fredholm pair.
\end{proof}  
 
\indent Next follows the preparation needed for the second characterization.\par

\newtheorem*{rem5.3}{Remark 5.3} 
\begin{rem5.3}\rm Let $X$ and $Y$ be two Banach spaces
and $S\in L(X,Y)$. Then, given another Banach space $Z$, it is possible 
to define the left and right multiplication
operators induced by $S$, that is
\begin{align*}
&L_S\colon L(Z,X)\to L(Z,Y),& &L_S(V)=SV,\\
&R_S\colon L(Y,Z)\to L(X,Z),&
&R_S(W)=WS,
\end{align*}
where $V\in L(Z,X)$ and $W\in L(Y,Z)$.\par

\indent It is clear that $\parallel L_S\parallel \le \parallel S\parallel$.
Furthermore, since $L_S(K(Z,X))\subseteq K(Z,Y)$ and 
$R_S(K(Y,Z))\subseteq K(X,Z)$, it is possible to introduce the operators
$$ 
\tilde{L}_S\colon C(Z,X)\to C(Z,Y),\hskip1cm \tilde{R}_S\colon C(Y,Z)\to C(X,Z),
$$
where 
\begin{align*}
&C(Z,X)=L(Z,X)/K(Z,X),& &C(Z,Y)=L(Z,Y)/K(Z,Y),\\
&C(Y,Z)=L(Y,Z)/K(Y,Z),& &C(X,Z)=L(X,Z)/K(X,Z),
\end{align*}
and the maps
$\tilde{L}_S$ and $\tilde{R}_S$ are the factorizations of $L_S$ and $R_S$
through the respective closed ideal of compact operators.\par 

\indent Similarly, if $T\in L(Y,X)$, then, given another Banach space $Z$, it is 
possible to define $L_T$ and $R_T$,
the left and right multiplication operators induced by $T$, that is
\begin{align*}
&L_T\colon L(Z,Y)\to L(Z,X),& &L_T(V)=TV,\\
&R_T\colon L(X,Z)\to L(Y,Z),&
&R_T(W)=WT,
\end{align*}
where $V\in L(Z,Y)$ and $W\in L(X,Z)$. Furthermore, as above, it is also 
possible to define
$$ 
\tilde{L}_T\colon C(Z,Y)\to C(Z,X),\hskip1cm \tilde{R}_T\colon C(X,Z)\to C(Y,Z),
$$
the factorizations of $L_T$ and $R_T$ through the respective closed ideal of 
compact operators. \par

\indent Next suppose that $R(ST)$ and $R(TS)$ are finite dimensional
subspaces of $Y$ and $X$ respectively. Then 
$$
\tilde{L}_S\tilde{L}_T=\tilde{L}_{ST}=0,\hskip1cm 
\tilde{L}_T\tilde{L}_S=\tilde{L}_{TS}=0,
$$
that is the pairs of operators $(\tilde{L}_S,\tilde{L}_T)$ and $(\tilde{L}_T,\tilde{L}_S)$
are chains, see \cite[10.3]{5} or \cite{6}.\par
\indent Finally, consider $U\in L(X_2,X_3)$ and $V\in L(X_1,X_2)$, where $X_1$, $X_2$ 
and $X_3$ are three Banach spaces, and  
suppose that $(U,V)$ is a chain, that is $UV=0$. The chain
$(U,V)$ is called \it{exact}\rm, if $R(V)=N(U)$. In addition, it is said that
$(U,V)$ is an \it{invertible chain}\rm, if there are continuous linear maps
$V_1\in L(X_2,X_1)$ and $U_1\in L(X_3,X_2)$ such that
$$
U_1U+VV_1=I,
$$
where $I$ denotes the identity map of $X_2$, see \cite[10.3.1]{5} or \cite{6}.\par  
\end{rem5.3}
\indent Next follows the second characterization of regular Fredholm pairs.\par

\newtheorem*{th5.4}{Theorem 5.4} 
\begin{th5.4}Let $X$ and $Y$ be two Banach spaces, and consider $S\in L(X,Y)$
and $T\in L(Y,X)$ two operators such that $R(ST)$ and $R(TS)$ are finite dimensional 
subspaces of $Y$ and $X$ respectively. With the same notations of Remark 5.3,
the following assertions are equivalent:\par
i) the pair $(S,T)$ belongs to $RP(X,Y)$;\par
ii) the operators $S$ and $T$ are regular, and 
$(\tilde{L}_S,\tilde{L}_T)$ $($resp. $(\tilde{L}_T,\tilde{L}_S)$$)$ is an invertible chain 
for any Banach space $Z$;\par
iii) the operators $S$ and $T$ are regular, and  
$(\tilde{L}_S,\tilde{L}_T)$ $($resp. $(\tilde{L}_T,\tilde{L}_S)$$)$ is an  invertible chain 
for the Banach space $X$ $($resp. $Y$$)$;\par
iv) the operators $S$ and $T$ are regular, and 
$(\tilde{L}_S,\tilde{L}_T)$ $($resp. $(\tilde{L}_T,\tilde{L}_S)$$)$ is an exact chain
for the Banach space $X$ $($resp. Y$)$ .\par
\indent Similarly, the following assertions are equivalent:\par
i) the pair $(S,T)$ belongs to $RP(X,Y)$;\par
ii) the operators $S$ and $T$ are regular, and $(\tilde{R}_S,\tilde{R}_T)$  $($resp. 
$(\tilde{R}_T,\tilde{R}_S)$$)$ is an invertible chain for any Banach space $Z$;\par
iii) the operators $S$ and $T$ are regular, and  
$(\tilde{R}_S,\tilde{R}_T)$ $($resp. $(\tilde{R}_T,\tilde{R}_S)$$)$ is an  invertible chain 
for the Banach space $Y$ $($resp. $X$$)$ ;\par
iv) the operators $S$ and $T$ are regular, and $(\tilde{R}_S,\tilde{R}_T)$ $($resp. 
$(\tilde{R}_T,\tilde{R}_S)$$)$ is an exact chain
for the Banach space $Y$ $($resp. X$)$ .\par 
\end{th5.4}
\begin{proof}
\indent First of all, observe that since $R(ST)$ and $R(TS)$
are finite dimensional Banach space, there exist $\mathcal{X}$ and
$\mathcal{Y}$, two closed subspaces of $X$ and $Y$ respectively, 
such that $X=\mathcal{X}\oplus R(TS)$
and $Y=\mathcal{Y}\oplus R(ST)$. Moreover, if $S$ and $T$ are presented as 
matrices, that is if
$$
S=\begin{pmatrix}
S_{11}&0\\
S_{21}&S_{22}
\end{pmatrix},
\hskip1cm
T=\begin{pmatrix}
T_{11}&0\\
T_{21}&T_{22}
\end{pmatrix},
$$
and if $X/R(TS)$ and $Y/R(ST)$ are identified with $\mathcal{X}$ and
$\mathcal{Y}$ respectively, then the maps $\overline{S}$ and $\overline{T}$
in Remark 5.1 can be identified with $S_{11}$ and $T_{11}$ respectively.
\par 
\indent In addition, if $Z$ is an arbitrary Banach space, since 
\begin{align*}
&L(Z,X)=L(Z,\mathcal{X})\oplus L(Z,R(TS)),& &L(Z,Y)=L(Z,\mathcal{Y})\oplus 
L(Z,R(ST)),\\
&K(Z,X)=K(Z,\mathcal{X})\oplus K(Z,R(TS)),& &K(Z,Y)=K(Z,\mathcal{Y})\oplus 
K(Z,R(ST)),
\end{align*}
then, 
$$
C(Z,X)=C(Z,\mathcal{X}),\hskip1cm C(Z,Y)=C(Z,\mathcal{Y}).
$$
\indent Furthermore, it is clear that 
$$
\tilde{L}_S=\tilde{L}_{S_{11}},\hskip1cm \tilde{L}_T=\tilde{L}_{T_{11}}.
$$ 

\indent Now well, if $(S,T)$ belongs to $RP(X,Y)$, then according to Theorem 5.2
and to the above identifications, $(S_{11},T_{11})$ is a regular Fredholm symmetrical 
pair. Therefore, according to \cite[10.6.2]{5},
there are operators $S_1$ and $S_2$ in $L(\mathcal{Y},\mathcal{X})$,
$T_1$ and $T_2$ in $L(\mathcal{X},\mathcal{Y})$, and two operators
with finite dimensional rank, $K_1\in L(\mathcal{X})$ and $K_2\in L(\mathcal{Y})$,
such that
$$
S_1S_{11}+T_{11}T_1=I_1-K_1,\hskip1cm T_2T_{11}+S_{11}S_2=I_2-K_2,
$$
where $I_1$ and $I_2$ denote the identity maps of $\mathcal{X}$ and 
$\mathcal{Y}$ respectively.\par
\indent Consequently,
$$
\tilde{L}_{S_1}\tilde{L}_{S_{11}} +\tilde{L}_{T_{11}}\tilde{L}_{T_1}=\mathcal{I}_1,
\hskip1cm \tilde{L}_{T_2}\tilde{L}_{T_{11}} +\tilde{L}_{S_{11}}\tilde{L}_{S_2}=
\mathcal{I}_2,
$$
where $\mathcal{I}_1$ and $\mathcal{I}_2$ denote the identity maps of 
$C(Z,\mathcal{X})$ and $C(Z,\mathcal{Y})$ respectively. \par
\indent However, since $C(Z,\mathcal{X})=C(Z,X)$, $C(Z,\mathcal{Y})=C(Z,Y)$,
$\tilde{L}_{S_{11}}=\tilde{L}_S$ and $\tilde{L}_{T_{11}}=\tilde{L}_T$,
then $(\tilde{L}_S,\tilde{L}_T)$ and $(\tilde{L}_T,\tilde{L}_S)$
are invertible chains in the sense of \cite[10.3.1]{5} or \cite{6}.\par

\indent It is clear that ii) implies iii) and iii) implies iv). \par

\indent Next suppose that iv) holds. Since $S$ and $T$ are regular maps, 
there are operators $S'\in L(Y,X)$ and $T'\in L(X,Y)$ such that
$S=SS'S$ and $T=TT'T$. \indent

\indent Now well, if, as above, $X$ and $Y$ are decomposed as direct sums
$$
X=\mathcal{X}\oplus R(TS),\hskip1cm Y=\mathcal{Y}\oplus R(ST),
$$
and if $S$ and $S'$ are presented as matrices, that is if
$$
S=\begin{pmatrix}
S_{11}&0\\
S_{21}&S_{22}
\end{pmatrix},
\hskip1cm
S'=\begin{pmatrix}
S_{11}'&S_{12}'\\
S_{21}'&S_{22}'
\end{pmatrix},
$$
then a straightforward calculation proves that
$$
S_{11}=S_{11}S_{11}'S_{11} + S_{11}S_1,
$$
where $S_1\in L(\mathcal{X})$ is an operator whose range is finite dimensional.\par 
\indent On the other hand, 
since $C(\mathcal{X})=C(\mathcal{X},\mathcal{X})=C(X,X)=C(X)$, $C(\mathcal{X},
\mathcal{Y})=C(X,Y)$,
$\tilde{L}_{S_{11}}=\tilde{L}_S$ and $\tilde{L}_{T_{11}}=\tilde{L}_T$,
the chain $(\tilde{L}_{S_{11}},\tilde{L}_{T_{11}})$ is exact.
However, if $\mathcal{I}$ denotes the identity of $C(\mathcal{X})$, since
$$
\tilde{L}_{S_{11}}(\mathcal{I}-[S_{11}'S_{11}])=0,
$$
then there is $B\in L(\mathcal{X},\mathcal{Y})$ such that
$$
\tilde{L}_{T_{11}}([B])= \mathcal{I}-[S_{11}'S_{11}],
$$
that is 
$$
T_{11} B+S_{11}'S_{11}=I-K,
$$
where $I$ is the identitiy map of $\mathcal{X}$, and $K\in L(\mathcal{X})$ is
a compact operator.\par
\indent Now well, since $(S_{11},T_{11})$ is a chain, that is $R(T_{11})\subseteq 
N(S_{11})$, and since 
$$
T_{11}B(N(S_{11}))=(I-K)(N(S_{11})),
$$
then $(I-K)(N(S_{11}))\subseteq N(S_{11})$. However, $N(S_{11})$ is a Banach space 
and $I-K$ is a Fredholm operator in $N(S_{11})$. Consequently, 
$\dim N(S_{11})/(I-K)(N(S_{11}))$ is finite dimensional, and since $(I-K)(N(S_{11}))
\subseteq R(T_{11})\subseteq N(S_{11})$, then $\dim N(S_{11})/R(T_{11})$ is finite.\par
\indent A similar argument proves that $\dim N(T_{11})/R(S_{11})$ is finite. 
Therefore, $(S_{11}, T_{11})$ is a Fredholm symmetrical pair.
However, according to the above identifications and to Remark 5.1, $(S,T)\in P(X,Y)$, 
and since $S$ and $T$ are regular operators, $(S,T)$ is a regular Fredholm pair.\par
\indent Similar arguments prove the second part of the theorem.\par
\end{proof}

\newtheorem*{th5.5}{Theorem 5.5}
\begin{th5.5} Let $X$ and $Y$ be two Banach spaces and consider $S\in L(X,Y)$
and $T\in L(Y,X)$, two regular operators such that $R(ST)$ and $R(TS)$ are finite
dimensional subspaces of $Y$ and $X$ respectively. Then, with the same notations
of Remark 5.3 and Theorem 5.4, if 
$S'$ and $T'$ are generalized inverses for $S$ and $T$ respectively, 
necessary and sufficient for $(S,T)$ to belong to $RP(X,Y)$ is that
$$
\tilde{L}_{S'}\tilde{L}_S + \tilde{L}_T\tilde{L}_{T'}\hbox{  and  }
\tilde{L}_{T'}\tilde{L}_T +\tilde{L}_S\tilde{L}_{S'}
$$
are invertible operators for any Banach space $Z$.\par
\indent Similarly, necessary and sufficient for $(S,T)$ to belong to $RP(X,Y)$ is that
$$
\tilde{R}_{S'}\tilde{R}_S + \tilde{R}_T\tilde{R}_{T'} \hbox{  and  }
\tilde{R}_{T'}\tilde{R}_T +\tilde{R}_S\tilde{R}_{S'}
$$
are invertible operators for any Banach space $Z$.\par
\end{th5.5}
\begin{proof}
\indent Since $R(ST)$ and $R(TS)$ are finite dimensional Banach spaces,
$(\tilde{L}_S,\tilde{L}_T)$ and $(\tilde{L}_T,\tilde{L}_S)$ are chains for 
any Banach space $Z$. Furthermore,
since $S=SS'S$ and $T=TT'T$, 
$$
\tilde{L}_S=\tilde{L}_S\tilde{L}_{S'}\tilde{L}_S,\hskip1cm
 \tilde{L}_T=\tilde{L}_T\tilde{L}_{T'}\tilde{L}_T,
$$
that is $\tilde{L}_{S'}$ and $\tilde{L}_{T'}$ are generalized inverses for $\tilde{L}_S$
and $\tilde{L}_T$ respectively.\par 

\indent Now well, acording to \cite[1]{6}, necessary and sufficient for
$$
\tilde{L}_{S'}\tilde{L}_S + \tilde{L}_T\tilde{L}_{T'}\hbox{ and  }
\tilde{L}_{T'}\tilde{L}_T +\tilde{L}_S\tilde{L}_{S'}
$$
to be invertible is the 
fact that $(\tilde{L}_S,\tilde{L}_T)$ and $(\tilde{L}_T,\tilde{L}_S)$ are invertible chains.
However, according to Theorem 5.4, this last assertion is equivalent to the fact
that $(S,T)$ belongs to $RP(X,Y)$.\par
\indent A similar argument proves the second part of the theorem.\par
\end{proof}

\vskip.3cm

Enrico Boasso\par
E-mail address: enrico\_odisseo@yahoo.it

\end{document}